\documentclass[twoside]{article}

\usepackage{PRIMEarxiv}

\usepackage{amssymb}

\usepackage{xcolor}
\usepackage{bm}
\usepackage{amsfonts} % For \mathbb
\usepackage[T1]{fontenc} % For the polish letter
\usepackage[caption=false]{subfig} % For subfigures
\usepackage[ruled,linesnumbered]{algorithm2e}
\SetKwInput{KwInput}{Input}        % Set the Input
\SetKwInput{KwOutput}{Output}      % set the Output
\usepackage{multicol}
\usepackage{url}
\usepackage{hyperref}
\usepackage{adjustbox}
\usepackage{caption}
\usepackage{rotating} 
\usepackage{amsthm}
\usepackage[style=numeric]{biblatex}
\usepackage{fontawesome}

\usepackage{tikz}
\usetikzlibrary{arrows.meta, decorations.markings}
\usepackage{tikz-network}
\usepackage{pgfplots}
\newcommand{\pasep}[1]{P_{\mbox{\tiny \textit{ASEP}}}^{#1}}

\newcommand{\erre}{\mathbb{R}}

\newcommand{\asep}[1]{\text{ASEP}(\bm{#1})}
\newcommand{\ssep}[1]{\text{SSEP}(\bm{#1})}
\newcommand{\atsp}[1]{\text{ATSP}(\bm{#1})}
\newcommand{\stsp}[1]{\text{STSP}(\bm{#1})}
\newcommand{\verasep}[1]{\mathcal{X}^{#1}_{\mbox{\tiny \textit{ASEP}}}}
\newcommand{\verintasep}[1]{\mathcal{T}^{#1}_{ASEP}}

\newcommand{\rev}{}
\newcommand{\revv}[1]{}

\newtheorem{thm}{Theorem}
\newtheorem{lemma}[thm]{Lemma}
\newtheorem{corollary}[thm]{Corollary}
\newtheorem{conj}{Conjecture}
\newtheorem{example}{Example}
\newtheorem{definition}{Definition}
\newtheorem{remark}{Remark}

\title{On the integrality Gap of Small Asymmetric Traveling Salesman Problems: A Polyhedral and Computational Approach}

\author{
  Eleonora Vercesi\textsuperscript{ \faEnvelopeO}, Luca Maria Gambardella \\
  Faculty of Informatics,\\
  Università della Svizzera italiana (USI)\\
  Lugano, Switzerland
  \texttt{\{eleonora.vercesi,luca.gambardella\}@usi.ch} \\
   \AND
  Janos Barta, Monaldo Mastrolilli \\
  Dipartimento Tecnologie innovative,\\
  Scuola universitaria professionale della Svizzera italiana (SUPSI)\\
  Lugano, Switzerland
  \texttt{\{janos.barta, monaldo.mastrolilli\}@supsi.ch} \\
  \AND
  Stefano Gualandi \\
  Department of Mathematics ``Felice Casorati''\\
  University of Pavia\\
  Pavia, Italy\\
  \texttt{stefano.gualandi@unipv.it} \\
}

\begin{document}
\maketitle

\begin{abstract}
In this paper, we investigate the integrality gap of the Asymmetric Traveling Salesman Problem (ATSP) with respect to the linear relaxation given by the Asymmetric Subtour Elimination Problem (ASEP) for \rev{instances with $n$ nodes, where $n$ is small}. In particular, we focus on the geometric properties and symmetries of the ASEP polytope ($\pasep{n}$)  and its vertices. The polytope's symmetries are exploited to design a heuristic pivoting algorithm \rev{to search} vertices where the integrality gap is maximized. 
Furthermore, a general procedure for the extension of vertices from $\pasep{n}$ to $\pasep{n + 1}$ is defined. The generated vertices improve the known lower bounds of the integrality gap for $ 16 \leq n \leq 22$ and, provide small hard-to-solve ATSP instances.
\end{abstract}


\section{Introduction}
\rev{The use of integer linear programming in modeling has been instrumental in both the practical solution of combinatorial optimization problems (See, e.g. ~\cite{applegate1998solution}) and the formulation and analysis of approximation algorithms (See, e.g ~\cite{vaziranibook})}. 
When solving an instance, an important concept is the \emph{integrality gap} \rev{(IG)} with respect to the linear relaxation, which is the maximum ratio between the solution quality of the Integer Linear Program (ILP) and its Linear Program (LP) relaxation.
For many integer linear programs,  the integrality gap of the linear relaxation is equal to the approximation ratio of the best algorithm, as well as the hardness of the approximation ratio, and essentially represents the inherent limits of the considered relaxation. To some extent, instances having a large integrality gap are the hard instances for the class of approaches based on linear programming. With this respect, the following facts can be observed from the literature. 

\begin{itemize}
    \item Proving integrality gaps for LP relaxations of NP-hard optimization problems is a difficult task that is usually undertaken on a case-by-case basis (e.g., see the case of the Vertex Cover \cite{singh2019integrality}).
    Very few and limited attempts have been made to use computer-assisted analysis for this difficult goal.
    \item For several notable and important examples, like the Traveling Salesman Problem (TSP), our integrality gap comprehension has resisted the persistent attack during the last decades of many researchers. For the Symmetric \rev{Traveling} Salesman Problem (STSP), several advances have been made in specific cases, such as the STSP having only costs 1 and 2 \cite{mnich2018improved}, and the cubic and subcubic STSP \cite{boyd2011tsp}, where the integrality gap is proved to be equal to $\frac{4}{3}$. For the asymmetric case, \rev{\cite{charikar2004integrality,elliott2008integrality}} \rev{exhibit two different families of ATSP instances of increasing size with an integrality gap asymptotically tending to 2. This, in turn, implies} that the lower bound for the integrality gap is 2, and no further improvements have been made. More specifically, \cite{elliott2008integrality} also provides specific lower bounds for $n \leq 15$, \rev{where $n$ is the number of nodes}. This is the first time an IG greater than $\frac{4}{3}$ is shown for $ n = 9$. In terms of the upper bound, the best-known bound is constant \rev{and} equal to 22 \cite{traub2020improved}.
\end{itemize}
Within this paper, we investigate the integrality gap 
\emph{computationally} for the Asymmetric Traveling Salesman Problem. 
More specifically, we investigate and exploit aspects of polyhedral theory to reduce the integrality gap search space.
Then, we combine this information to effectively search the pruned space.
This approach allows us to get computer-aided construction of bad instances, namely, problem instances having large integrality gap values, for the considered $LP$ relaxation, which improves upon the best-known lower bounds for small $n$ (see Charikar, Goemans, and Karloff  \cite{charikar2004integrality} and Elliot-Magwood~\cite{elliott2008integrality}). 
The aim is to obtain a new and better understanding of lower bounds on the corresponding integrality gaps. \rev{As a side result}, we present hard-to-solve ATSP instances\rev{.}

\paragraph{The problem} 
Given a weighted directed graph, a Hamiltonian cycle is a cycle that visits each \rev{node} exactly once.
The Asymmetric Traveling Salesman Problem (ATSP) involves finding a directed Hamiltonian cycle of minimum cost.
In practice, the currently most efficient exact algorithm for solving ATSP is based on an Integer Linear Programming (ILP) \rev{formulation} \rev{\cite{fischetti2007exact,fischetti1997polyhedral,roberti2012models}}.
Let $K_n = (V, A)$ be the complete directed graph with $n$ nodes and $m = n(n-1)$ arcs, that is, $V:=\{1,\dots,n\}$ and $A:=\{(i,j) \mid i,j \in V, i\neq j\}$, each having a weight $c_{ij} \in \erre^+$.
Whenever $c_{ij} = c_{ji}$ for all arcs $(i,j)$, then we have a Symmetric TSP instance.
Whenever $c_{ij} \neq c_{ji}$, but the cost vector satisfies the triangle inequality $c_{ij} \leq c_{ik} + c_{kj}, \forall i,j,k \in V$, we have a pseudo-quasi metric TSP instance, since the cost vector induces a {\it pseudo-quasi metric} \cite{lawvere1973metric}.
Note that any instance of ATSP on a complete graph $K_n$ is completely defined by its cost vector $\bm{c} \in \mathbb{R}^m_+$.
Given $S \subset V$, let $\delta(S)$ be the cut induced by $S$, namely $\delta(S) := \{ (i, j) \mid i \in S, j \not \in S\}$. For convenience, we denote $\delta(i) = \delta(\{i\})$.

In this paper, we study the LP relaxation of the Dantzig-Fulkerson-Johnson (DFJ) formulation \cite{dantzig1954solution} for solving the ATSP for small values of $n$.
\begin{align}
\min \quad & \sum_{(i,j) \in A} c_{ij} x_{ij} \label{obj} \\
\text{s.t.} \quad & \sum_{i \in V} x_{ij} = 1 && \forall j \in V \label{indegree}\\
 \quad & \sum_{j \in V} x_{ij} = 1 && \forall i \in V \label{outdegree}\\
\quad & \sum_{\rev{(i, j)}\in \delta(S)} x_{ij} \geq 1 && \forall S \subset V \mbox{ such that } 2 \leq \vert S  \vert \leq  n - 2 \label{subtour}\\
\quad & x_{ij} \in \{0, 1\} && \forall (i,j) \in A, \label{int}
\end{align}
where $x_{ij} \in \{0,1\}$ is equal to 1 if arc $(i,j)$ \rev{belongs to an} optimal cycle, and 0 otherwise.
Constraints \eqref{indegree}--\eqref{outdegree} are the in-degree and out-degree constraints, which force each node to have exactly one predecessor and one successor. Constraints \eqref{subtour} are the \emph{Subtour Elimination Constraints}, which avoid the presence of sub-cycles (subtours). For a survey on ATSP formulations, we refer the reader to \cite{roberti2012models}.
In the literature, the solution of the LP relaxation of \eqref{obj}--\eqref{int} is called the \emph{Asymmetric Subtour Elimination Problem} (ASEP), while the corresponding feasibility region is the \emph{ASEP Polytope}, which is defined as follows. 
\begin{equation}
    \label{pasep}
    \pasep{n} := \{ \bm{x} \in \mathbb{R}^{m} \; \vert \;  \text{\eqref{indegree}--\eqref{subtour}}, \bm{x} \geq 0\}. 
\end{equation}

The Integrality Gap for the ATSP on $n$ nodes is
\begin{equation}
    \alpha_n := \sup_{c \in \mathbb{R}^{m}_+} \frac{\text{ATSP}(\bm{c})}{\text{ASEP}(\bm{c})},
\end{equation}
 where $\atsp{c}$ is the optimal value of \eqref{obj}--\eqref{int} and $\asep{c}$ is the optimal value of its LP relaxation.
In general, the IG for the ATSP is 
\[\alpha:= \sup_{n \in \mathbb{N}} \alpha_n.\]

In \cite{elliott2008integrality}, it is shown that for the general ATSP, a non-negative cost vector $\bm{c}$ exists such that $ \atsp{c} > 0$ and $\asep{c} = 0$.
This implies that the integrality gap tends to $+\infty$.
For this reason, the works studying the integrality gap of ATSP always assume that the cost vector satisfies the triangle inequality $c_{ij} \leq c_{ik} + c_{kj}, \forall i,j,k \in V$ and, hence, the costs induce a pseudo-quasi metric.
Herein, we restrict our attention to pseudo-quasi-metric ATSP, but we will call them ATSP for short, to be consistent with the notation of the literature.

A long-standing conjecture stated that $\alpha \leq \frac{4}{3}$ (e.g., see \cite{carr2004held}). 
However, this conjecture was disproved \rev{in  \cite{charikar2004integrality,elliott2008integrality}, where it has been shown that} $\alpha \geq 2$ by presenting two families of ATSP instances with $\alpha(\bm{c}) \to 2$ for $n \to \infty$. 
For $n \leq 3$, we have that $\alpha_n=1$, while for $4 \leq n \leq 7$ the exact value was computed in \cite{elliott2008integrality}.
The authors of \rev{\cite{BoydE05,elliott2008integrality}} provide also lower bounds of $\alpha_n$ for $n \leq 15$.
In particular, they show that $\alpha > \frac{4}{3}$, since they prove (compute) that $\alpha_{n}\geq \frac{11}{8}$ for $n=9$.
Lower bounds were also provided by Charikar et al. \cite{charikar2004integrality} for arbitrarily large $n$, but they are rather weak for $n\leq 25$. 
We remark that the literature is more extensive for the Symmetric TSP \rev{\cite{benoit2008finding,HOUGARDY2014495,hougardy2020hard,vercesi2023generation,zhong2021lower}}.
For instance, \rev{\cite{hougardy2020hard,vercesi2023generation,zhong2021lower}} propose STSP instances that have a large integrality gap and are hard to solve for the state-of-the-art \rev{STSP} solver Concorde \cite{applegate1998solution}.
However, no theoretical proof is available for the exact value of the integrality gap, neither for the STSP nor the ATSP.
Furthermore, the relation between computational complexity and large integrality gap is still unclear even in the symmetric case.

\paragraph{Main contributions.} This paper has three main contributions. 
First, we identify, for the first time, a group of symmetries of the ASEP polytope.
We exploit this symmetry to design pivoting strategies that explore vertices of $\pasep{n}$.
Second, we provide new lower bounds for the integrality gap $\alpha_n$ for the ATSP for $16 \leq n \leq 22$ by using a new pivoting algorithm which exploits the symmetries of the vertices of $\pasep{n}$, combined with an inductive algorithm that generates vertices of $\pasep{n+1}$ from a vertex of $\pasep{n}$. 
Our bounds improve those provided in \rev{\cite{charikar2004integrality} and \cite{elliott2008integrality}}. 
Third, by using our new inductive algorithm, we generate hard ATSP instances, where complexity is measured with respect to the Concorde solver for STSP~\cite{applegate1998solution} after an appropriate transformation of the ATSP instance, \rev{and with respect to two dedicated ATSP solvers \cite{fischetti2007exact,fischetti1992additive}}.

The outline of this paper is as follows.
Section \ref{sec:background} reviews the background material. 
Section \ref{sec:sym} studies the symmetries of the polytope $\pasep{n}$.
In Section \ref{sec:four}, we explain how we use the algebraic structure introduced to perform a symmetry-breaking heuristic pivoting.
Furthermore, we introduce a new operator that generates vertices of $\pasep{n + 1}$ from a vertex of $\pasep{n}$.
In Section \ref{sec:results} we present the results of our approach, by analyzing the structure of the obtained vertices and exhibiting new lower bounds for $16 \leq n \leq 22$.
Finally, in Section \ref{sec:final} we conclude the paper with a perspective on future works.

%%%%%%%%%%%%%%%%%%%%%%%%%%%%%%
% Background
%%%%%%%%%%%%%%%%%%%%%%%%%%%%%%
\section{Background material}
\label{sec:background}
A key subproblem of our approach is the computation of the maximum integrality gap over all possible pseudo-quasi metric cost vectors $\bm{c} \in \mathbb{R}^{m}_+$ for a single vertex of the ASEP polytope. This subproblem was first introduced in \cite{benoit2008finding} for the \rev{STSP} and later in \cite{elliott2008integrality} for the \rev{ATSP}. The main idea is to divide the cost vector $\bm{c}$ by the optimal value $\text{ATSP}(\bm{c})$, \rev{obtaining $\bm{c'}$ such that $\text{ATSP}(\bm{c}') = 1$. This transformation \rev{leaves} unchanged any tour leading to an optimal solution.  Hence, one can consider without loss of generality the definition of IG as to $\alpha_{n}:=\sup _{\bm{c}\text{ s.t  ATSP}(\bm{c}) = 1} \frac{1}{\text{ASEP}(\bm{c})}$, which is equivalent to solve $\inf_{\bm{c}\text{ s.t  ATSP}(\bm{c}) = 1} \text{ASEP}(\bm{c})$}. Note that the operation of dividing the costs by $\text{ATSP}(\bm{c})$ maintains the triangle inequalities and preserves the value of the IG.

The problem of computing the maximum $\alpha_{n}$ of $K_{n}$ for the pseudo-quasi metric (pq-metric) ATSP is as follows \cite{elliott2008integrality}.

\begin{equation}
  \frac{1}{\alpha_{n}}:=\min _{\substack{\bm{c} \text { is pq-metric, } \\ \text{ATSP}(\bm{c})=1}} \text{ASEP}(\bm{c})=\min _{\bm{x} \in \pasep{n}} \min _{\substack{\text { is pq-metric, } \\ \text{ATSP}(\bm{c})=1}} \bm{x}^{T} \bm{c}. \label{eq:gap_problem}  
\end{equation}

To solve the problem \eqref{eq:gap_problem} for a fixed (small) value of $n$, the authors in \cite{elliott2008integrality} enumerate the vertices of $\pasep{n}$, and for each vertex $\bm{x} \in \pasep{n}$, they solve $\min \left\{\bm{c}^{T} \bm{x} \mid\right.$ $\bm{c}$ is pq-metric, $\text{ATSP}(\bm{c})=1\}$. In \cite{elliott2008integrality}, the latter inner problem is called $\text{Gap}(\bm{x})$. Intuitively, \rev{$\text{Gap}(\bm{x})$ allows us to determine, among all cost vectors that have an optimal ASEP solution at point $\bm{x}$, the one that results in the largest integrality gap. Since the number of vertices of $\pasep{n}$ is finite, we only need to explore a fixed number of cases. With a slight abuse of notation, we refer to ``the integrality gap of a vertex'' as the integrality gap induced by the cost vector that maximizes the integrality gap among those with an optimal ASEP solution at $\bm{x}$.}

\rev{The definition of $\text{Gap}(\bm{x})$ relies heavily on the slackness compatibility conditions, which ensure that an optimal primal-dual pair satisfies complementary slackness.}

\rev{For this reason, the definition of $\text{Gap}(\bm{x})$ involves the dual variables of the linear formulation \eqref{indegree}--\eqref{subtour}, where $\bm{x} \geq 0$. Let $y_{j}^{\text {in }}$ and $y_{i}^{\text {out }}$ be the dual variables associated with the in-degree constraint \eqref{indegree} and out-degree constraint \eqref{outdegree}, respectively. Let $d_S$ be the dual variables associated with the subtour elimination constraint \eqref{subtour}. Finally, let $c_{ij}$ denote the cost variable on arc $(i, j)$. The function $\text{Gap}(\bm{x})$ is defined as below. For further details, the interested reader can check \cite{elliott2008integrality}.}

\begin{align}
     \text{Gap}(\bm{x}):=&\min \sum_{(i, j) \in A} x_{i j} c_{i j} \label{gap:obj_func}\\
\text { s.t. }&  \sum_{(i, j) \in A} z_{i j} c_{i j} \geq 1 \quad \forall \bm{z} \in \mathcal{T}_{A S E P}^{n} \label{constr:atsp_opt}\\
& c_{i j} \leq c_{i k}+c_{j k} \quad \forall i, j, k \in V \\
& c_{i j} \geq 0 \quad \forall(i, j) \in A \\
& c_{i j}-y_{i}^{\text {out }}-y_{j}^{\text {in }}-\sum_{S \in \mathcal{S}_{i j}(\bm{x})} d_{S} \geq 0 \quad \forall(i, j) \in A \backslash A(\bm{x}) \label{const:gap_1}\\
& c_{i j}-y_{i}^{\text {out }}-y_{j}^{\text {in }}-\sum_{S \in \mathcal{S}_{i j}\left(\bm{x}^{*}\right)} d_{S}=0 \quad \forall(i, j) \in A(\bm{x}) \label{const:gap_2}\\
& d_{S} \geq 0 \quad \forall S \in \bigcup_{(i, j) \in A} \mathcal{S}_{i j}(\bm{x}), \label{const:gap_3}
\end{align}
where $\mathcal{T}_{\text {ASEP }}^{n}$ is the collection of every possible Hamiltonian cycle of $K_{n}$, and $\bm{z} \in\{0,1\}^{m}$ are incidence vectors of elements of $\mathcal{T}_{\text {ASEP }}^{n}$. We remark that the variables are $c_{i j}$, \rev{$y_i^{\text{in}}, y_i^{\text{out}}, d_S$}, while \rev{$\bm{x}, \bm{z}$} are given. Constraints \eqref{constr:atsp_opt} ensure that
the optimal solution $\bm{c}^{*}$ of $\text{Gap}(\bm{x})$ yields $\text{ATSP}\left(\bm{c}^{*}\right)=1$. Constraints \eqref{const:gap_1}--\eqref{const:gap_3} derive from dualizing the linear program relaxation of \eqref{indegree}--\eqref{subtour}. Here, we define the following subset of nodes $S_{i j(\bm{x})}$ and of arcs $A(\bm{x})$ as follows: $\mathcal{S}_{i j}(\bm{x})=\{S \subset A|(i, j) \in \delta(S), \bm{x}(\delta(S))=1,2 \leq| S \mid \leq n-2\}$ and $A(\bm{x})=\left\{(i, j) \in A \mid x_{i j}>0\right\}$. Constraints \eqref{const:gap_1}--\eqref{const:gap_3} ensure that

$$
\bm{x}^{*} \in \arg \min _{\bm{x} \in \pasep{n}} \bm{c}^{*} \bm{x} \Leftrightarrow \bm{c}^{*} \in \arg \min _{\bm{c}, \bm{y}, \bm{d} \text { satisfy \eqref{constr:atsp_opt}--\eqref{const:gap_3}}} \text{Gap}\left(\bm{x}^{*}\right) .
$$

Thus, the $\arg \min \bm{c}^{*}$ of the program $\text{Gap}\left(\bm{x}^{*}\right)$ is such that, once $\text{ASEP}\left(\bm{c}^{*}\right)$ is solved, the minimum is attained precisely at $\bm{x}^{*}$ (for details, see \cite{elliott2008integrality}).

%%%%%%%%%%%%%%%%%%%%%%%%%%%%%%
%Symmetric group
%%%%%%%%%%%%%%%%%%%%%%%%%%%%%%
\section{Symmetry group of $\pasep{n}$}
\label{sec:sym}
Polyhedral aspects of the TSP have been extensively studied in the past decades. In particular, in \rev{\cite{balinski1974on,grotschelbook,padberg_rao}} the polyhedral properties of the convex hull of the integer vertices, the so-called ``natural polytope'', of the STSP and the ATSP have been investigated starting from the assignment polytope. 
More specifically, it has been shown that the diameter of the assignment polytope is two, which in turn implies that the diameter of the TSP polytope is at most 2.
\rev{However, since the valid formulation for the TSP examined in this paper -- namely, the one proposed by Dantzig, Fulkerson, and Johnson \cite{dantzig1954solution} -- is not integral, our focus is on the relationship between optimizing only on the integer points and optimizing over $\pasep{n}$ (linear relaxation).}
In this section, we focus in particular on the fractional vertices and the symmetry properties of the subtour elimination polytope $\pasep{n}$.

As remarked in \cite{elliott2008integrality}, the vertices of $\pasep{n}$ can be subdivided into equivalence classes through permutations of the nodes of $K_n$. In this section, we define explicitly a symmetry group of $\pasep{n}$ based on permutations. The classes of isomorphic vertices turn out to be the orbits generated by this symmetry group.

By observing the definition of $\pasep{n}$, an important consideration can be made:
Due to the structure of the constraints, they are not affected by a permutation of the node indices. More precisely, it is well known that any relabeling of the nodes $i \in V$ induces just an internal permutation of the constraint groups (\ref{indegree}), (\ref{outdegree}), and (\ref{subtour}), which leaves the feasible region unchanged.
Consequently, it is possible to identify a group of symmetries of the polytope $\pasep{n}$ induced by the symmetric group $S_n$ of permutations of the nodes $i\in V$.
To describe the symmetries of $\pasep{n}$ explicitly, it is useful to convert the feasible solutions into matrix form. More precisely, we can rewrite any feasible solution $\bm{x}\in \pasep{n}$ as a matrix $\bm{X} \in [0,1]^{n\times n}$, with the corresponding components $x_{ij}$ for $i\neq j$ and setting $x_{ij}=0$ for $i=j$.

Now, consider the symmetric group of permutations on $n$ elements $S_n$.
Let $\pi \in S_n$, such that $\pi : i \mapsto \pi(i)$, be a permutation of the nodes $i\in V$. For any feasible solution $\bm{x}\in\pasep{n}$, with matrix representation $\bm{X}$, it is possible to generate a new feasible solution $\bm{x}'$, expressed by a matrix $\bm{X'}$, by applying the permutation $\pi$ to the nodes of the graph $K_n$. It must hold
\begin{equation} \label{permutation_x}
    x'_{\pi(i)\pi(j)} = x_{ij} \quad \forall i,j\in V.
\end{equation}
This means that $\bm{X'}$ is obtained by permuting rows and columns of $\bm{X}$ according to $\pi$. \rev{Given a permutation of $n$ elements $\pi\in S_n$, } let us recall the definition of \emph{permutation matrix} \rev{as the matrix $\bm{P}_\pi = (p_{ij})$ such that}
\[
p_{ij} = \begin{cases}
1 & \mbox{if } i=\pi(j) \\
0 & \mbox{otherwise}.
\end{cases}
\]
In the product $\bm{P_\pi X}$, the permutation matrix $\bm{P_\pi}$ permutes the rows of $\bm{X}$ according to $\pi$.
Since also the columns of $\bm{X}$ have to be permuted, we apply $\bm{P_\pi}$ to the transpose of $\bm{P_\pi X}$.
Thus
\[ \bm{X}' = \left(\bm{P_\pi}(\bm{P_\pi}\bm{X})^T\right)^T
= (\bm{P_\pi}\bm{X})\bm{P_\pi}^T { = \bm{P_\pi}\bm{X}\bm{P_\pi}^T.}
\]  
A notion that we will widely use is the \emph{isomorphism between digraphs.}
Hence, let us recall some helpful definitions.

\begin{definition}[Support digraph] 
Let  $\bm{x} \in \pasep{n}$. The weighted support digraph $\bm{x}$ is the graph $H(\bm{x})$ defined on the set of nodes $V$, having an arc $(i,j)$ with the weight $x_{ij}$, if and only if $x_{ij}>0$. 
\end{definition}

\begin{definition}[Isomorphism between weighted graphs]
Two graphs  $D=(V,A_D)$ and $F=(V,A_F)$ on $n$ nodes are isomorphic if there exists a permutation $\pi$ of $V$ such that $(u,v)\in A_D \iff (\pi(u),\pi(v)) \in A_F$. Furthermore, the weights \rev{$\bm{c}^D, \bm{c}^F$} of the edges must satisfy \rev{$c_{\pi(u)\pi(v)}^F=c_{uv}^D, \forall u,v \in V$}.
\label{def:graph_isom}
\end{definition}

\begin{example}
Consider $\bm{x}\in \pasep{4}$ defined by 
\[ \bm{x} = \left(\frac{1}{2}, 0,\frac{1}{2},\frac{1}{2},\frac{1}{2}, 0,\frac{1}{2}, 0,\frac{1}{2}, 0,\frac{1}{2},\frac{1}{2} \right)^T.\]
Its matrix version is hence 
\[ \bm{X} = \left(\begin{array}{cccc}0 & \frac{1}{2} & 0 & \frac{1}{2} \\ \frac{1}{2} & 0 & \frac{1}{2} & 0 \\ \frac{1}{2} & 0 & 0 & \frac{1}{2} \\ 0 & \frac{1}{2} & \frac{1}{2} & 0\end{array}\right). \]
This feasible solution is associated with the support graph in Figure \ref{fig:esempio_4}, left.
Let $\pi = (0 \; 1 \; 2\;  3)$, that is the permutation such that $\pi(0) = 1, \; \pi(1) = 2, \; \pi(2) = 3, \;  \pi(3) = 0.$
Thus, in this case 
\[\bm{P}_\pi=\left(\begin{array}{llll}0 & 0 & 0 & 1 \\ 1 & 0 & 0 & 0 \\ 0 & 1 & 0 & 0 \\ 0 & 0 
& 1 & 0\end{array}\right)\]
and we  obtain
\[\bm{X}^{\prime}=\bm{P}_\pi \bm{X} \bm{P}_\pi^T=\left(\begin{array}{cccc}0 & 0 & \frac{1}{2} & \frac{1}{2} \\ \frac{1}{2} & 0 & \frac{1}{2} & 0 \\ 0 & \frac{1}{2} & 0 & \frac{1}{2} \\ \frac{1}{2} & \frac{1}{2} & 0 & 0\end{array}\right),\]
which corresponds to the graph in Figure \ref{fig:esempio_4} on the right. 
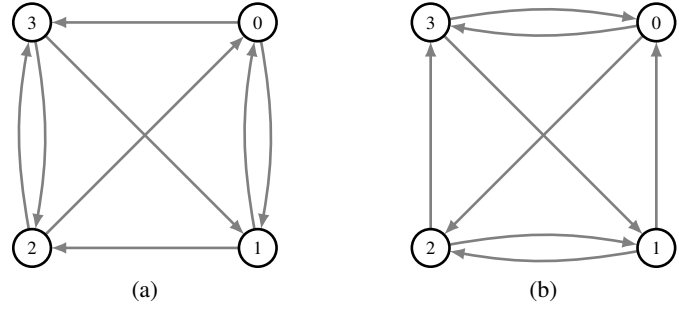
\begin{figure*}[!t]
\centering
   \subfloat[]{\label{rev}
      \begin{tikzpicture}
\Vertex[x=4,y=4,,size=0.5,color=white,opacity=1,fontcolor=black,label=0]{1}
\Vertex[x=4,y=1,,size=0.5,color=white,opacity=1,fontcolor=black,label=1]{2}
\Vertex[x=1,y=1,,size=0.5,color=white,opacity=1,fontcolor=black,label=2]{3}
\Vertex[x=1,y=4,,size=0.5,color=white,opacity=1,fontcolor=black,label=3]{4}
\Edge[,lw=1.0,color={127.5,127.5,127.5},Direct,RGB, bend=10](1)(2)
\Edge[,lw=1.0,color={127.5,127.5,127.5},Direct,RGB](1)(4)
\Edge[,lw=1.0,color={127.5,127.5,127.5},Direct,RGB, bend=10](2)(1)
\Edge[,lw=1.0,color={127.5,127.5,127.5},Direct,RGB](4)(2)
\Edge[,lw=1.0,color={127.5,127.5,127.5},Direct,RGB, bend=10](3)(4)
\Edge[,lw=1.0,color={127.5,127.5,127.5},Direct,RGB](3)(1)
\Edge[,lw=1.0,color={127.5,127.5,127.5},Direct,RGB](2)(3)
\Edge[,lw=1.0,color={127.5,127.5,127.5},Direct,RGB, bend=10](4)(3)
\end{tikzpicture}}
    \hspace{1.5cm}
   \subfloat[]{\label{rev_sol}
      \begin{tikzpicture}
\Vertex[x=4,y=4,size=0.5,color=white,opacity=1,fontcolor=black,label=0]{1}
\Vertex[x=4,y=1,size=0.5,color=white,opacity=1,fontcolor=black,label=1]{2}
\Vertex[x=1,y=1,size=0.5,color=white,opacity=1,fontcolor=black,label=2]{3}
\Vertex[x=1,y=4,,size=0.5,color=white,opacity=1,fontcolor=black,label=3]{4}
\Edge[,lw=1.0,color={127.5,127.5,127.5},Direct,RGB, bend=10](4)(1)
\Edge[,lw=1.0,color={127.5,127.5,127.5},Direct,RGB](2)(1)
\Edge[,lw=1.0,color={127.5,127.5,127.5},Direct,RGB, bend=10](1)(4)
\Edge[,lw=1.0,color={127.5,127.5,127.5},Direct,RGB](1)(3)
\Edge[,lw=1.0,color={127.5,127.5,127.5},Direct,RGB, bend=10](3)(2)
\Edge[,lw=1.0,color={127.5,127.5,127.5},Direct,RGB](4)(2)
\Edge[,lw=1.0,color={127.5,127.5,127.5},Direct,RGB](3)(4)
\Edge[,lw=1.0,color={127.5,127.5,127.5},Direct,RGB, bend=10](2)(3)
\end{tikzpicture}}
   \caption{Two isomorphic vertices obtained via vertex permutation. Each arc is weighted $\frac{1}{2}$.}\label{fig:esempio_4}
\end{figure*}
\end{example}

For any node relabeling $\pi \in S_n$ we can define the corresponding permutation of a solution $\bm{X}$ as

\begin{align} \label{g_pi}
    g_\pi: \; & [0,1]^{n\times n} \to [0,1]^{n\times n} \\
    & \bm{X} \mapsto \bm{P}_\pi \bm{X} \bm{P}_\pi^T. \nonumber
\end{align} 

With a slight abuse of notation, we will use both $g_\pi(\bm{x})$ and $g_\pi(\bm{X})$ interchangeably, denoting the \emph{row and column permutation} of $\bm{X}$, according to $\pi$.

First of all, we observe that the map is well defined in $\pasep{n}$, namely $\bm{x}\in \pasep{n}$ implies $g_\pi(\bm{x}) \in \pasep{n}$.
Hence, all the constraints are satisfied also for $g_\pi(\bm{x})$ by a simple ``shuffle'' of the rows: degree constraints $i$ become degree constraints $\pi(i)$ and subtour elimination constraints associated to $\delta(S), \; S = \{s_1, \ldots, s_r\}$ become $\delta(\pi(S)), \; \pi(S) := \{\pi(s_1), \ldots, \pi(s_r)\}$.
The isomorphism of vertices, observed in \cite{elliott2008integrality}, can now be extended to the whole polytope $\pasep{n}$. Hence, it trivially follows
\begin{lemma} \label{lemma_iso}
    Let $x\in \pasep{n}$ and $\pi \in S_n$. The support graphs $H(\bm{x})$ and $H(g_\pi(\bm{x}))$ are isomorphic.
\end{lemma}
Let \rev{$G(n) = \{g_\pi \; \vert \; \pi \in S_n  \}$ be} the set of all transformations $g_\pi$. The next theorem shows that $G(n)$ is a group of symmetries of the polytope.

\begin{thm} \label{thm_symmetry_group}
    $G(n)$ is a group of isometries acting on of $\pasep{n}$. 
\end{thm}
\begin{proof}
{We begin by showing the closure of $G(n)$ under   composition. Let ${\pi_1}, {\pi_2} \in S_n$. By equation (\ref{g_pi}) we have
\begin{eqnarray} \label{G(n)_closure}
    g_{\pi_1} g_{\pi_2}(\bm{X}) &=&g_{\pi_1}\left(\bm{P}_{\pi_2} \bm{X} \bm{P}_{\pi_2}^{T}\right) = \bm{P}_{\pi_1} \bm{P}_{\pi_2} \bm{X} \bm{P}_{\pi_2}^{T} \bm{P}_{\pi_1}^{T} \\
    &=& \left(\bm{P}_{\pi_1 \pi_2} \right) \bm{X} \left(\bm{P}_{\pi_1 \pi_2} \right)^{T}
    = g_{\pi_1 \pi_2}(\bm{X}) \in G(n). \nonumber
\end{eqnarray}
It is not difficult to verify that the symmetric group $S_n$ induces the group structure on $G(n)$. In particular, by equation (\ref{G(n)_closure}) we have that the identity element of $G(n)$ is $g_{id}$ and the inverse element of $g_\pi$ is $g_{\pi^{-1}}$. Moreover, equation (\ref{G(n)_closure}) states that $G(n)$ and $S_n$ are isomorphic groups.
By equation (\ref{permutation_x}) it can be observed that for any feasible solution $\bm{x} \in \pasep{n}$ the solution $\bm{x}'=g_\pi(\bm{x}) \in \pasep{n}$ is obtained by a permutation of the components of $\bm{x}$ based on $\pi$. Therefore, it is immediately clear that $g_\pi$ preserves the Euclidean distance, that is  $\vert \vert g_\pi(\bm{x}) - g_\pi(\bm{y}) \vert \vert = \vert \vert \bm{x} - \bm{y} \vert \vert \ , \forall \bm{x},\bm{y} \in \pasep{n}, \pi \in S_n$. }
\end{proof}

As explained in Section \ref{sec:background}, our main goal is the computation of the integrality gap of vertices. The following corollary focuses on the action of $G(n)$ on the set of vertices, which will be denoted throughout the manuscript with $\verasep{n}$. 

\begin{corollary}
Let $\bm{x} \in \verasep{n}$ and $g_\pi \in G(n)$. Then  $g_\pi(\bm{x}) \in \verasep{n}.$
\label{thm:orbits}
\end{corollary}
\begin{proof} 
{$g_\pi$ is an isometry and isometries map vertices into vertices.}
\end{proof}
Since the group $G(n)$ acts on the set of vertices $\verasep{n}$, it is of interest to study the $\emph{orbit}$ of each vertex $\bm{x} \in \verasep{n}$, that is the set
\begin{equation} \label{orb}
O_{\bm{x}} = \{ g_\pi(\bm{x}) \; \vert \; g_\pi \in G(n) \},
\end{equation}
and the \emph{stabilizer} of $G(n)$ with respect to $\bm{x} \in \verasep{n}$, that is the subgroup of $G(n)$ defined as
\[ G_{\bm{x}} = \{g_\pi \in G(n) \; \vert \; g_\pi(\bm{x}) = \bm{x}\}.\]

{Combining Lemma \ref{lemma_iso} and \eqref{orb}, we can conclude that vertices of the polytope $\pasep{n}$ belonging to the same orbit have isomorphic support graphs. In other terms, the isomorphism classes of vertices introduced in \cite{elliott2008integrality} are the orbits of the vertices with respect to the group $G(n)$.}

By the so-called \emph{orbit-stabilizer theorem} we have the relation
$\vert G(n) \vert = \vert G_{\bm{x}} \vert \vert O_{\bm{x}} \vert$.
Further details on the orbit-stabilizer theorem can be found, for instance, in \cite{artin2011algebra}.
A natural question concerns the role of the stabilizer of each vertex, and whether it leads to some implications in combinatorial questions, such as the integrality gap. For example, what could be the relationship between the stabilizer of a vertex and the maximum integrality gap achievable at that vertex? The case of the integer vertices is particularly interesting. As already pointed out, integer vertices of $\pasep{n}$ {correspond to Hamiltonian cycles} of $K_n$. However, Hamiltonian cycles differ from each other only by a relabeling of the nodes. Therefore, we can prove the following property.

\begin{lemma} \label{lem_int_vertices}
Let $\bm{x}\in \verintasep{n}$. Then, it holds $O_{\bm{x}} = \verintasep{n}$.
%\ele{Non noto}
\end{lemma}
\begin{proof}
Let $\bm{x_1},\bm{x_2} \in \verintasep{n}$. Since $\bm{x_1}$ and $\bm{x_2}$ {correspond to Hamiltonian cycles in $K_n$}, there exists a permutation $\pi$, such that $g_\pi(\bm{x_1})=\bm{x_2}$. It follows that $\bm{x_1}$ and $\bm{x_2}$ belong to the same orbit of $G(n)$.
\end{proof}

A  consequence of Lemma \ref{lem_int_vertices} is that integer vertices build a unique orbit $O_{\bm{x}}$ with $|O_{\bm{x}}|=(n-1)!$. By the orbit-stabilizer theorem it follows that $|G_{\bm{x}}| = n, \ \forall \bm{x}\in \verintasep{n}$.

More specifically, it holds the following.

\begin{lemma}
    If $\bm{x} \in \verintasep{n}$, then $G_{\bm{x}}$ is the group
    \[\langle (1 \; 2 \ldots n) \rangle \]
    of all the cyclic permutations that are
    \[\tau_k(i) = (i + k) \mod n.\]
\end{lemma}
\begin{proof}
    As these permutation\rev{s} are the $k$-shift of nodes, $k\in \{0, 1, 2, \ldots, n\}$, it holds 
    \[ \tau_{k_1} \circ \tau_{k_2} = \tau_{k_3} \qquad k_3 = {(k_1 + k_2) \mod{n}}.\]
\end{proof}

Unfortunately, as further discussed in Section \ref{sec:results_symmetries}, the stabilizer of a fractional vertex $\bm{x}$ of $\pasep{n}$ does not seem to follow a general recipe. 

In Section \ref{sec:final}, we conjecture that vertices having a large integrality gap are the ones having large but not trivial stabilizers; hence, it could be relevant to know in advance the structure of the stabilizer to guess the ``promising'' vertices.

It is worth remarking that the orbits of $G(n)$ form large classes of isomorphic vertices; hence, $\pasep{n}$ can be considered a highly symmetric polytope.

The intrinsic equivalence of vertices belonging to the same orbit is visible in the following lemma.
\begin{lemma}
Let $\bm{x},\bm{y} \in \verasep{n}$ such that $\bm{y}\in O_{\bm{x}}$, then $Gap(\bm{x}) = Gap(\bm{y})$. 
\end{lemma}
\begin{proof}
If $\bm{x}, \bm{y} \in O_{\bm{x}}$, then there exists $g_{\pi} \in G(n)$, such that $\bm{x} = g_\pi(\bm{y})$. 
Let $\bm{c}^* \in \arg \min Gap(\bm{x})$.
\[Gap(\bm{x}) = \bm{c}^{*T} \bm{x} = g_{\pi}(\bm{c}^{*})^T g_{\pi}(\bm{x}) = Gap(\bm{y}).\]
 The first equation holds by definition, while the second equation remains unchanged even if both terms of the scalar product are permuted. The final equation is derived from the fact that if there exists a $\bm{c}'$ value such that $\bm{c}^{'^T}\bm{x}' < \bm{c}^{*^T} \bm{x}$, then $g_{\pi}^{-1}(\bm{c}')$ would provide a solution of lower cost than $\bm{c}^{*^T}$, thus rendering the latter non-optimal.
\end{proof}

Note that this result has already been proved in \cite{elliott2008integrality}, using a different strategy.
Finally, we add a result that justifies our symmetry-breaking simplex algorithm presented in the next section.

\begin{definition}
    The set of all vertices adjacent to $\bm{x} \in \verasep{n}$ is called the {\bf neighborhood of} $\bm{x}$ and is denoted by $\mathcal{N}(\bm{x})$.
\end{definition}

\begin{lemma}\label{lemma:neigh}
    Let  $\bm{x},\bm{y} \in \verasep{n}$, $\pi\in S_n$. Then, $\bm{y}\in \mathcal{N}(\bm{x}) \Leftrightarrow g_\pi(\bm{y}) \in  \mathcal{N}(g_\pi(\bm{x}))$. 
\end{lemma}
\begin{proof}
First, we use the transformation $\pi$ to sort the rows of the constraint matrix, by mapping each degree constraint associated with $i$ to $\pi(i)$, and each subtour elimination constraint associated with $S = \{s_1, \ldots, s_f\}$ to $\pi(S):=  \{s_1, \ldots, s_f\}$ and each non-negative constraint accordingly. With a slight abuse of notation, we will call $\pi(i)$ the mapping of the $i^{th}$ constraint.
Note that $\bm{y}$ shares exactly $r(n) - 1$ linearly independent and tight constraints with $\bm{x}$, let  $i_1, \ldots, i_{n_y - 1}, i_{n_y}$ are this set of constraints of  $\bm{y}$ and $i_1, \ldots, i_{n_y - 1}, i_{n_x}$ is the one of  $\bm{x}^3$, then, $\pi\left(i_1\right), \ldots, \pi\left(i_{n_y - 1}\right), \pi\left(i_{n_y}\right)$ is associated to $g_\pi(\bm{y}$ and $\pi\left(i_1\right), \ldots, \pi\left(i_{n_y-1}\right), \pi\left(i_{n_x}\right)$ is hence associated with $g_{\pi}(\bm{x})$. Thus, $g_\pi(\bm{y}) \in \mathcal{N}(g_{\pi})(\bm{x}))$.
\end{proof}
{Substantially, Lemma \ref{lemma:neigh} states that the maps $g_\pi$ preserve the adjacency of vertices. In other words, vertices belonging to the same orbit are equivalent also in terms of their neighborhoods.}

%%%%%%%%%%%%%%%%%%%%%%%%%%%%%%%%%%
% New vertices
%%%%%%%%%%%%%%%%%%%%%%%%%%%%%%%%%%
\section{Computing vertices with \rev{a large} integrality gap}\label{sec:four}

The orbits of the vertices of $\pasep{n}$ introduced in the previous section are here used to design a computational strategy for heuristically generating vertices of $\pasep{n}$.
In the next paragraphs, first, we introduce our pivoting algorithm that exploits the vertex symmetries to avoid the exploration of isomorphic vertices.
Then, we introduce a new iterative procedure that generates vertices of $\pasep{n+1}$ starting from vertices of $\pasep{n}$.

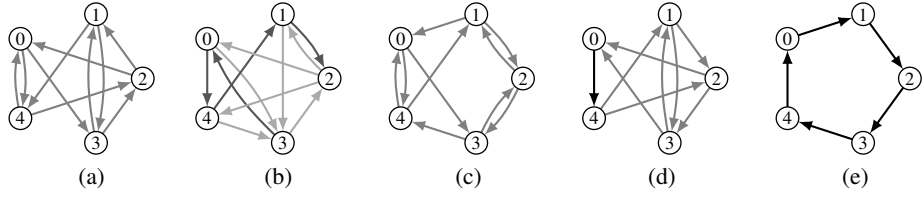
\begin{figure*}[!t]
\centering 
\subfloat[]{\begin{tikzpicture}[scale=0.30]
\clip (-0.4,-0.4) rectangle (7,7);
\renewcommand*{\VertexLineWidth}{0.5pt}

\Vertex[x=0.289,y=4.761,size=0.3,color=white,opacity=1,label=0,fontcolor=black]{0}
\Vertex[x=5.711,y=3.000,size=0.3,color=white,opacity=1,label=2,fontcolor=black]{2}
\Vertex[x=3.640,y=5.850,size=0.3,color=white,opacity=1,label=1,fontcolor=black]{1}
\Vertex[x=0.289,y=1.239,size=0.3,color=white,opacity=1,label=4,fontcolor=black]{4}
\Vertex[x=3.640,y=0.150,size=0.3,color=white,opacity=1,label=3,fontcolor=black]{3}
\Edge[,lw=0.8,color={127.5,127.5,127.5},Direct,RGB](0)(3)
\Edge[,lw=0.8,color={127.5,127.5,127.5},Direct,RGB, bend=10](0)(4)
\Edge[,lw=0.8,color={127.5,127.5,127.5},Direct,RGB, bend=10](3)(1)
\Edge[,lw=0.8,color={127.5,127.5,127.5},Direct,RGB](3)(2)
\Edge[,lw=0.8,color={127.5,127.5,127.5},Direct,RGB, bend=10](4)(0)
\Edge[,lw=0.8,color={127.5,127.5,127.5},Direct,RGB](4)(2)
\Edge[,lw=0.8,color={127.5,127.5,127.5},Direct,RGB, bend=10](1)(3)
\Edge[,lw=0.8,color={127.5,127.5,127.5},Direct,RGB](1)(4)
\Edge[,lw=0.8,color={127.5,127.5,127.5},Direct,RGB](2)(0)
\Edge[,lw=0.8,color={127.5,127.5,127.5},Direct,RGB](2)(1)
\end{tikzpicture}}%     without .tex extension
~    
\subfloat[]{\begin{tikzpicture}[scale=0.30]
\clip (-0.4,-0.4) rectangle (7,7);
\renewcommand*{\VertexLineWidth}{0.5pt}
\Vertex[x=0.289,y=4.761,size=0.3,color=white,opacity=1,label=0,fontcolor=black]{0}
\Vertex[x=5.711,y=3.000,size=0.3,color=white,opacity=1,label=2,fontcolor=black]{2}
\Vertex[x=3.640,y=5.850,size=0.3,color=white,opacity=1,label=1,fontcolor=black]{1}
\Vertex[x=0.289,y=1.239,size=0.3,color=white,opacity=1,label=4,fontcolor=black]{4}
\Vertex[x=3.640,y=0.150,size=0.3,color=white,opacity=1,label=3,fontcolor=black]{3}
\Edge[,lw=0.8,color={170.0,170.0,170.0},Direct,RGB, bend=10](0)(3)
\Edge[,lw=0.8,color={85.0,85.0,85.0},Direct,RGB](0)(4)
\Edge[,lw=0.8,color={85.0,85.0,85.0},Direct,RGB, bend=10](3)(0)
\Edge[,lw=0.8,color={170.0,170.0,170.0},Direct,RGB](3)(2)
\Edge[,lw=0.8,color={85.0,85.0,85.0},Direct,RGB](4)(1)
\Edge[,lw=0.8,color={170.0,170.0,170.0},Direct,RGB](4)(3)
\Edge[,lw=0.8,color={85.0,85.0,85.0},Direct,RGB, bend=10](1)(2)
\Edge[,lw=0.8,color={170.0,170.0,170.0},Direct,RGB](1)(3)
\Edge[,lw=0.8,color={170.0,170.0,170.0},Direct,RGB](2)(0)
\Edge[,lw=0.8,color={170.0,170.0,170.0},Direct,RGB, bend=10](2)(1)
\Edge[,lw=0.8,color={170.0,170.0,170.0},Direct,RGB](2)(4)
\end{tikzpicture}%     without .tex extension
}
~
\subfloat[]{
        \begin{tikzpicture}[scale=0.30]
\clip (-0.4,-0.4) rectangle (7,7);
\renewcommand*{\VertexLineWidth}{0.5pt}

\Vertex[x=0.289,y=4.761,size=0.3,color=white,opacity=1,label=0,fontcolor=black]{0}
\Vertex[x=5.711,y=3.000,size=0.3,color=white,opacity=1,label=2,fontcolor=black]{2}
\Vertex[x=3.640,y=5.850,size=0.3,color=white,opacity=1,label=1,fontcolor=black]{1}
\Vertex[x=0.289,y=1.239,size=0.3,color=white,opacity=1,label=4,fontcolor=black]{4}
\Vertex[x=3.640,y=0.150,size=0.3,color=white,opacity=1,label=3,fontcolor=black]{3}
\Edge[,lw=0.8,color={127.5,127.5,127.5},Direct,RGB](0)(3)
\Edge[,lw=0.8,color={127.5,127.5,127.5},Direct,RGB, bend=10](0)(4)
%\Edge[,lw=0.8,color={127.5,127.5,127.5},Direct,RGB, bend=10](3)(0)
\Edge[,lw=0.8,color={127.5,127.5,127.5},Direct,RGB](4)(1)
\Edge[,lw=0.8,color={127.5,127.5,127.5},Direct,RGB, bend=10](4)(0)
\Edge[,lw=0.8,color={127.5,127.5,127.5},Direct,RGB](3)(4)
\Edge[,lw=0.8,color={127.5,127.5,127.5},Direct,RGB, bend=10](1)(2)
\Edge[,lw=0.8,color={127.5,127.5,127.5},Direct,RGB](1)(0)
\Edge[,lw=0.8,color={127.5,127.5,127.5},Direct,RGB, bend=10](2)(1)
\Edge[,lw=0.8,color={127.5,127.5,127.5},Direct,RGB, bend=10](2)(3)
\Edge[,lw=0.8,color={127.5,127.5,127.5},Direct,RGB, bend=10](3)(2)
\end{tikzpicture}%     without .tex extension
}~
\subfloat[]{
        \begin{tikzpicture}[scale=0.30]
\clip (-0.4,-0.4) rectangle (7,7);
\renewcommand*{\VertexLineWidth}{0.5pt}

\Vertex[x=0.289,y=4.761,size=0.3,color=white,opacity=1,label=0,fontcolor=black]{0}
\Vertex[x=5.711,y=3.000,size=0.3,color=white,opacity=1,label=2,fontcolor=black]{2}
\Vertex[x=3.640,y=5.850,size=0.3,color=white,opacity=1,label=1,fontcolor=black]{1}
\Vertex[x=0.289,y=1.239,size=0.3,color=white,opacity=1,label=4,fontcolor=black]{4}
\Vertex[x=3.640,y=0.150,size=0.3,color=white,opacity=1,label=3,fontcolor=black]{3}
\Edge[,lw=0.8,color={0.0,0.0,0.0},Direct,RGB](0)(4)
\Edge[,lw=0.8,color={127.5,127.5,127.5},Direct,RGB](4)(1)
\Edge[,lw=0.8,color={127.5,127.5,127.5},Direct,RGB](4)(2)
\Edge[,lw=0.8,color={127.5,127.5,127.5},Direct,RGB](1)(2)
\Edge[,lw=0.8,color={127.5,127.5,127.5},Direct,RGB, bend=10](1)(3)
\Edge[,lw=0.8,color={127.5,127.5,127.5},Direct,RGB](2)(0)
\Edge[,lw=0.8,color={127.5,127.5,127.5},Direct,RGB](2)(3)
\Edge[,lw=0.8,color={127.5,127.5,127.5},Direct,RGB](3)(0)
\Edge[,lw=0.8,color={127.5,127.5,127.5},Direct,RGB, bend=10](3)(1)
\end{tikzpicture}%     without .tex extension
}~
\subfloat[]{
        \begin{tikzpicture}[scale=0.30]
\clip (-0.4,-0.4) rectangle (7,7);
\renewcommand*{\VertexLineWidth}{0.5pt}
\Vertex[x=0.289,y=4.761,size=0.3,color=white,opacity=1,label=0,fontcolor=black]{0}
\Vertex[x=5.711,y=3.000,size=0.3,color=white,opacity=1,label=2,fontcolor=black]{2}
\Vertex[x=3.640,y=5.850,size=0.3,color=white,opacity=1,label=1,fontcolor=black]{1}
\Vertex[x=0.289,y=1.239,size=0.3,color=white,opacity=1,label=4,fontcolor=black]{4}
\Vertex[x=3.640,y=0.150,size=0.3,color=white,opacity=1,label=3,fontcolor=black]{3}
\Edge[,lw=0.8,color={0.0,0.0,0.0},Direct,RGB](0)(1)
\Edge[,lw=0.8,color={0.0,0.0,0.0},Direct,RGB](1)(2)
\Edge[,lw=0.8,color={0.0,0.0,0.0},Direct,RGB](2)(3)
\Edge[,lw=0.8,color={0.0,0.0,0.0},Direct,RGB](3)(4)
\Edge[,lw=0.8,color={0.0,0.0,0.0},Direct,RGB](4)(0)
\end{tikzpicture}%     without .tex extension
}
\caption{Support graphs of the 5 isomorphism classes for $n=5$. The arcs have a grey level depending on the value of the corresponding $x_i$. See Table \ref{tab:orbit_5} for details: while (a) and (c) have all $x_i=\frac{1}{2}$ (light grey); in the support graph (b) 4 arcs correspond to $x_i=\frac{2}{3}$ (dark grey) and 7 to $x_i=\frac{1}{3}$; (d) 8 arcs have $x_i=\frac{1}{2}$ and 1 arc $x_i=1$; (e) all $x_i=1$.}
\label{fig:n5repr}
\end{figure*} 

\subsection{Pivoting by symmetry-breaking}
\label{sec:pivoting}
In this subsection, we illustrate the new pivoting algorithm, which attempts to avoid the exploration of new vertices that are isomorphic to vertices already visited. The pivoting algorithm will be denoted by $\text{Pivoting}(\bm{x}, T)$, where the meaning of the variable $T$ will be clarified later.

The main idea is simple: we start with a basic feasible solution, that is a vertex $\bm{x} \in \pasep{n}$, and we explore all vertices in the neighborhood $\mathcal{N}(\bm{x})$ one at a time, by enumerating (or by sampling) the possible pivoting steps for that vertex.
If the new vertex obtained by pivoting is not isomorphic to any vertex already explored, then we solve the optimization problem \eqref{gap:obj_func}--\eqref{const:gap_3}
to find the maximal integrality gap for that new vertex, and we record the corresponding orbit.
We iterate the neighborhood search either with a complete enumeration for small values of $n$ (e.g., $n \leq 8)$, or with a random sampling strategy for $n > 8$.
Our procedure is iterative, namely, it continues to iterate vertex-by-vertex, exploring each time the neighborhood of each vertex.
The input parameters are:
\begin{itemize}
\item $M$, the maximum number of iterations of the algorithm, equivalent to the maximum number of vertices we pivot on.
\item $T_{tot}$, timelimit for the whole iterations.
\item  $T_{it}$, timelimit for a single iteration.
\end{itemize}
\rev{Parameter} $T_{it}$ balance the tradeoff between \emph{exploration} and \emph{exploitation}. 
Small values of $T_{it}$ allow for exploring only a few elements adjacent to a given vertex $\bm{x}$ and quickly moving on to the next neighboring vertex to be explored. On the other hand, high values of $T_{it}$ insist heavily on the neighborhood of a vertex $\bm{x}$.
We continue to iterate the procedure until either the timelimit $T_{tot}$ is hit or the maximum number $M$ of iteration\rev{s} is reached.

We named this algorithm ``explore/exploit'' for two reasons: the parameter $T_{it}$ governs the exploitation of a single vertex. For small values of $T_{it}$, we do not focus much on the neighborhood of a single vertex, instead preferring to pivot on multiple vertices. However, for large values of $T_{it}$, we continue to build a single $\mathcal{N}(\bm{x})$, having less time to explore different areas of $\pasep{n}$.

Hence, given a time limit of $T$,  our function 
\begin{equation}
\mathcal{N}'(\bm{x}) = \text{Pivoting}(\bm{x}, T)
\label{eq:piiv_func}
\end{equation}
returns only a subset of $\mathcal{N}(\bm{x})$, namely the adjacent vertices that the strategy is capable of finding within a time limit.

\subsection{Generating vertices using loop breaking procedure}
\label{sec:lambda-loops}
In this subsection, we present our new iterative algorithm to compute a vertex of $\pasep{n+1}$ by starting from a vertex of $\pasep{n}$.
First, we recover two definitions and a lemma from \cite{elliott2008integrality} that we use to prove our new result.
Then, we introduce our loop-breaking procedure.
\begin{definition}[\cite{elliott2008integrality}, Chap. 3]
\label{def:ts}
    Let $S \subset V$ and $\bm{x} \in \pasep{n}$. If $\delta(\bm{x}(S)) :=\sum_{(i, j) \in \delta(S)} x_e = 1$, then, $S$ is called a \textbf{tight set}.
\end{definition}
{Note that a tight set is a subset of vertices that induces a cut that satisfies a subtour elimination constraint with equality.}
\begin{definition}[\cite{elliott2008integrality}, Chap. 3]
\label{def:bb_collapse}
    Let $S \subset V$ be a tight set and $\bm{x} \in \pasep{n + \vert S \vert - 1}$. Then, we can collapse the set $S$ into node $w$ as follows: 
\[ 
    (\bm{x} \downarrow_w (S))_{uv} = 
    \begin{cases} 
        \sum_{s \in S} x_{us} & \text{if } v = w, \\
        \sum_{s \in S} x_{sv} & \text{if } u = w, \\
        x_{uv} & \text{otherwise.}
    \end{cases}
 \]
\end{definition}

\begin{lemma}[\cite{elliott2008integrality}, Prop. 3.3.1]
\label{lemma:prop331}
    Let $\bm{x} \in \pasep{n}$ and $S \subset V$ be a tight set of $\bm{x}$. Then, $\bm{x} \downarrow_w (S)$ belongs to $\pasep{n - \vert S \vert + 1}$. 
\end{lemma}
We are now ready to formally introduce $\lambda$-loops.
\begin{definition}[$\lambda$-loop]
Let $\bm{x} \in \pasep{n}$ and let $v_1, v_2 \in V$ such that $x_{v_1 v_2} = \lambda$ and $x_{v_2 v_1} = 1-\lambda$, then, we say that $\bm{x}$ \textbf{contains a $\lambda$-loop} $\overleftrightarrow{v_1 v_2}$.
\end{definition}

The following trivially follows.
\begin{lemma}
        Consider $\bm{x}$ a point in $\pasep{n}$ containing the $\lambda$-loop $\overleftrightarrow{v_1 v_2}$, then $S = \{v_1, v_2\}$ is a tight set.
\end{lemma}
\begin{proof}
    The amount of flow entering $v_1$ is $\lambda$, while the amount of flow entering $v_2$ is $1 - \lambda$. Hence, $\bm{x}(\delta(S)) = 1$.
\end{proof}

\subsubsection{A new $\lambda$-loop breaking procedure}
Our  algorithm for generating vertices of $\pasep{n + 1}$ from $\pasep{n}$ exploits $\lambda$-loop in the following $\lambda$-loop breaking procedure.

\begin{definition}\label{def:break_loops}
    Let $\bm{x} \in \pasep{n}$ that contains a $\lambda$-loop $\overleftrightarrow{v_1 v_2}$. Then, the \textbf{$\lambda$-loop breaking procedure} generate\rev{s} a point in $\erre^{(n+1)n}$ adding the node $v_3$ as follows :
    \[ 
    \left(\bm{x}  \uparrow^{v_3} (v_1, v_2)\right)_{uv} = 
    \begin{cases} 
        \lambda & \text{if } uv \in \{ v_1 v_3 , \;  v_3 v_2\}, \\
        1 - \lambda & \text{if } uv \in \{v_3 v_1 , \; v_2 v_3 \}, \\
        0 & \text{if } uv \in \{ v_1 v_2, \; v_2 v_1 \} \\
        0 & \text{if } v = v_3 \text{ and } u \not \in \{v_1, v_2\} \\
        x_{uv} & \text{otherwise.}
    \end{cases}
    \]
\end{definition}

%%%%% Figura
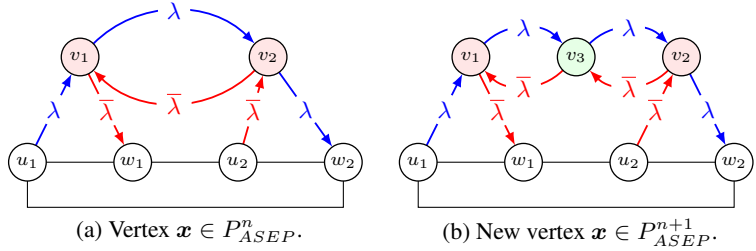
\begin{figure}[t!]
\centering
\subfloat[Vertex $\bm{x} \in P_{ASEP}^n$.]{
 \begin{tikzpicture}
\draw[draw=black] (0,0) rectangle (4.2,0.6);
\renewcommand*{\VertexLineWidth}{0.5pt}

\Vertex[x=0,y=0.6,size=0.5,color=white,opacity=1,label=$u_1$,fontcolor=black]{u1}
\Vertex[x=1.4,y=0.6,size=0.5,color=white,opacity=1,label=$w_1$,fontcolor=black]{w1}
\Vertex[x=2.8,y=0.6,size=0.5,color=white,opacity=1,label=$u_2$,fontcolor=black]{u2}
\Vertex[x=4.2,y=0.6,size=0.5,color=white,opacity=1,label=$w_2$,fontcolor=black]{w2}
\Vertex[x=0.7,y=2,size=0.5,color=red,opacity=0.1,label=$v_1$,fontcolor=black]{v1} \Vertex[x=3.2,y=2,size=0.5,color=red,opacity=0.1,label=$v_2$,fontcolor=black]{v2}    

\Edge[,lw=0.7,color=red,Direct,label={\small$\;\overline{\lambda}\;$}](u2)(v2)
\Edge[,lw=0.7,color=blue,Direct,label={\small$\;{\lambda}\;$}](v2)(w2)
\Edge[,lw=0.7,color=blue,Direct,label={\small$\;{\lambda}\;$}](u1)(v1)
\Edge[,lw=0.7,color=red,Direct,label={\small$\;\overline{\lambda}\;$}](v1)(w1)
\Edge[,lw=0.7,color=blue,Direct, bend=45,label={\small$\;\lambda\;$}](v1)(v2)
\Edge[,lw=0.7,color=red,Direct, bend=45,label={\small$\;\overline{\lambda}\;$}](v2)(v1)
\end{tikzpicture}
\label{fig:ext_left}
}~
\subfloat[New vertex $\bm{x} \in P_{ASEP}^{n+1}$.]{
 \begin{tikzpicture}
\draw[draw=black] (0,0) rectangle (4.2,0.6);
\renewcommand*{\VertexLineWidth}{0.5pt}

\Vertex[x=0,y=0.6,size=0.5,color=white,opacity=1,label=$u_1$,fontcolor=black]{u1}
\Vertex[x=1.4,y=0.6,size=0.5,color=white,opacity=1,label=$w_1$,fontcolor=black]{w1}
\Vertex[x=2.8,y=0.6,size=0.5,color=white,opacity=1,label=$u_2$,fontcolor=black]{u2}
\Vertex[x=4.2,y=0.6,size=0.5,color=white,opacity=1,label=$w_2$,fontcolor=black]{w2}
\Vertex[x=0.7,y=2,size=0.5,color=red,opacity=0.1,label=$v_1$,fontcolor=black]{v1}
\Vertex[x=3.5,y=2,size=0.5,color=red,opacity=0.1,label=$v_2$,fontcolor=black]{v2}
\Vertex[x=2.1,y=2,size=0.5,color=green,opacity=0.1,label=$v_3$,fontcolor=black]{v3} 
        
    \Edge[,lw=0.7,color=red,Direct,label={\small$\;\overline{\lambda}\;$}](u2)(v2)
    \Edge[,lw=0.7,color=blue,Direct,label={\small$\;{\lambda}\;$}](v2)(w2)
    \Edge[,lw=0.7,color=blue,Direct,label={\small$\;{\lambda}\;$}](u1)(v1)
    \Edge[,lw=0.7,color=red,Direct,label={\small$\;\overline{\lambda}\;$}](v1)(w1)
    
    \Edge[,lw=0.7,color=blue,Direct,bend=45,label={\small$\;{\lambda}\;$}](v1)(v3)
    \Edge[,lw=0.7,color=red,Direct,bend=45,label={\small$\;\overline{\lambda}\;$}](v3)(v1)
    \Edge[,lw=0.7,color=red,Direct,bend=45,label={\small$\;\overline{\lambda}\;$}](v2)(v3)
    \Edge[,lw=0.7,color=blue,Direct,bend=45,label={\small$\;{\lambda}\;$}](v3)(v2)
\end{tikzpicture}

%\draw[draw=black] (0,0) rectangle (3,0.6);
%\Vertex[x=0,y=0.6,size=0.5,color=white,opacity=1,label=$w_2$,fontcolor=black]{w2}
%\Vertex[x=1,y=0.6,size=0.5,color=white,opacity=1,label=$u_2$,fontcolor=black]{u2}
%\Vertex[x=2,y=0.6,size=0.5,color=white,opacity=1,label=$w_1$,fontcolor=black]{w1}
%\Vertex[x=3,y=0.6,size=0.5,color=white,opacity=1,label=$u_1$,fontcolor=black]{u1}
%\Vertex[x=0.7,y=2,size=0.5,color=white,opacity=1,label=$v_2$,fontcolor=black]{v2} \Vertex[x=2.3,y=2,size=0.5,color=white,opacity=1,label=$v_1$,fontcolor=black]{v1}    

%\Edge[,lw=0.7,color=red,Direct,label={\small$\overline{\lambda}$}](u2)(v2)
%\Edge[,lw=0.7,color=blue,Direct,label={\small${\lambda}$}](v2)(w2)
%\Edge[,lw=0.7,color=blue,Direct,label={\small${\lambda}$}](u1)(v1)
%\Edge[,lw=0.7,color=red,Direct,label={\small$\overline{\lambda}$}](v1)(w1)
%\Edge[,lw=0.7,color=blue,Direct, bend=45,label={\small$\lambda$}](v1)(v2)
%\Edge[,lw=0.7,color=red,Direct, bend=45,label={\small$\overline{\lambda}$}](v2)(v1)
\label{fig:ext_right}
}
\caption{Example of $\lambda$-loop breaking $(\bm{x} \uparrow^{v_3} (v_1, v_2))_{uv}$, with $\overline{\lambda} = 1-\lambda$.}\label{fig:extend}
\end{figure}

The $\lambda$-loop breaking procedure has the following properties.
\begin{lemma}
    $S = \{ v_1, v_3 \}$ is a tight set for $\bm{x} \uparrow^{v_3} (v_1, v_2)$.
\end{lemma}
\begin{proof}
Let $\bm{x}' = \bm{x}  \uparrow^{v_3} (v_1, v_2) $
\[\bm{x}^{\prime}(\delta(S))=\sum_{\substack{b \in \rev{V} \\ b \neq v_3}} x_{v_1 b}+\sum_{\substack{b \in \rev{V} \\ b \neq v_1}} x_{v_3 b} \]
given that
\[
\begin{aligned}
& \sum_{b \in V} x_{v_1 b} =1 \Rightarrow \sum_{\substack{b \in V \\
b \neq v_3}} x_{v_1 b} =1-x_{v_1 v_3}=1-\lambda, \\
& \sum_{b \in V} x_{v_3 b} =1 \Rightarrow \sum_{\substack{b \in V \\
b \neq v_1}} x_{v_3 b}=1-x_{v_3 v_1}=1-(1-\lambda)=\lambda, \\
& \Rightarrow \bm{x}^{\prime}(\delta(S))=1-\lambda+\lambda=1. \\
\end{aligned}\]
\end{proof}

\begin{lemma}\label{lemma:break_loops}
    If $\bm{x}$ is a vertex of $\pasep{n}$ that contains a $\lambda$-loop $\overleftrightarrow{v_1 v_2}$, then $\bm{x}  \uparrow^{v_3} (v_1, v_2)$ is a vertex of $\pasep{n + 1}$.
\end{lemma}
\begin{proof}
    Let $\bm{x}' := \bm{x}  \uparrow^{v_3} (v_1, v_2)$. Clearly $\bm{x}'$ is feasible. 
    To show that \rev{$\bm{x}$} is \rev{also} a vertex, assume by contradiction that there exists $\bm{y}', \bm{z}' \in \pasep{n + 1}$ and $\mu \in (0, 1)$ such that
    $\bm{x}' = \mu \bm{y}' + (1 - \mu)\bm{z}'$.
    Let $S = \{v_1, v_3\}$. Note that 
    \begin{equation}\label{eq:convex_comb}
        \bm{x}'(\delta(S)) = 1 = \mu \bm{y}'(\delta(S)) + (1 - \mu)\bm{z}'(\delta(S)).
    \end{equation}
    As $\bm{y}', \bm{z}'$ are feasible, we have $\bm{y}'(\delta(S)) \geq 1$ and $\bm{z}'(\delta(S)) \geq 1$. As the equality should hold, we have then  $\bm{z}'(\delta(S)) = \bm{y}'(\delta(S)) = 1$.
    Thus, we have
    \[ 
        \bm{y} := \downarrow_{v_1} \bm{y}' (S) \quad \text{ and } \quad \bm{z} := \downarrow_{v_1} \bm{z}' (S). 
    \]
    Thanks to Lemma \ref{lemma:prop331}, $\bm{z}, \bm{y} \in \pasep{n}$.
    Clearly, $\bm{x} = \mu \bm{y} + (1-\mu)\bm{z}$.
    This can be verified entry by entry. As an example, consider the case $e = uv_1$:
    \begin{eqnarray*}
        \mu y_{u v_1}+(1-\mu) z_{u v_1}&=&
  \mu\left(y_{u v_1}^{\prime}+y_{u v_3}^{\prime}\right)+(1-\mu)\left(z_{u v_1}^{\prime}+z_{u v_3}^{\prime}\right)  \\ 
        &=& \left[\mu y_{u v_1}^{\prime}+(1-\mu) z_{u v_1}^{\prime}\right]+\left[\mu y_{u v_3}^{\prime}+(1-\mu) z_{u v_3}^{\prime}\right] \\
&=& x_{u v_1}^{\prime}+x_{u v_3}^{\prime}=x_{u v_1}.
    \end{eqnarray*}
    Thus, $\bm{x}$ is not a vertex of $\pasep{n}$, and we get a contradiction. 
\end{proof}

\begin{remark}
Note that the converse unfortunately does not hold, e.g., if you collapse a $\lambda$-loop into a single \rev{node}, you are not guaranteed to obtain a vertex. 
This will be particularly relevant in Section \ref{sec:new_lb} where we collapse instead $\lambda$-loops: after doing that, we have to check that what we obtain is again a vertex.
\rev{However, it always leads to a feasible point thanks to Lemma \ref{lemma:prop331}.}
\end{remark}

%%%%%

\subsubsection{How we used this procedure}
\label{sec:how_to_use_lambda}
If we are able to compute an initial vertex \rev{$\bm{x}_0 \in \pasep{n}$} having a $\lambda-$loop, then, we can iteratively apply the $\lambda$-loop breaking procedure to obtain a sequence of vertices \rev{$\bm{x}_1, \bm{x}_{2}, \dots, \bm{x}_{t}$} that belongs to $\pasep{n+1}, \pasep{n+2}, \dots, \pasep{n+t}$, respectively.
Note that a similar idea was explored in \cite{elliott2008integrality}, where the author introduced a \emph{2-jack insertion} procedure, where a node satisfying precise hypotheses is replaced with a $\lambda$-loop, with the strict condition $\lambda = \frac{1}{2}$. 
Our strategy is more general since we allow any value of $\lambda \in (0, 1)$.
Our experimental results show that our procedure is very effective in finding vertices with large integrality gaps (see Section \ref{sec:impact_ll}).

\subsection{The full algorithm}
For a formal description of the proposed procedure, see Algorithm \ref{alg:expl}. The Pivoting function mentioned in line 9 of Algorithm \ref{alg:expl} is the one of Equation \eqref{eq:piiv_func}.
The whole procedure can be described in words as follows. \rev{First, let $T_{it}$ denote the time limit for exploring the neighborhood of a single vertex, $T_{tot}$ the overall time limit for the algorithm, and $M$ the maximum number of iterations allowed.}
Starting from $n = 5$ and the full collection of non-isomorphic vertices of $\pasep{5}$ (See Figure \ref{fig:n5repr}) that we can exhaustively generate using the software Polymake \cite{assarf2017computing}:

\begin{enumerate} 
\item We apply the procedure of $\lambda$-loop break obtaining vertices of $\pasep{n+1}$.
\item We complete each of these vertices with the slack variables.
\item We order all the vertices found in this way by the number of zeros, from the one with the fewest zeros to the one with the most zeros.
\item Starting from the first one, we begin to apply the Pivoting strategy as described in Section \ref{sec:pivoting} and Equation \ref{eq:piiv_func}. Specifically, we enumerate all possible combinations of variables that can form a feasible basis and attempt to include each of the nonbasic variables in the basis.
\item When we reach the time limit $T_{it}$, we have a subset of adjacent vertices denoted as $\mathcal{N}'(\bm{x}) \subseteq \mathcal{N}(\bm{x})$. At this point, two things can occur:
\begin{enumerate}
\item $\mathcal{N}'(\bm{x}) = \emptyset$: In this case, we move on to the next vertex in the ordered list and start again from step 4.
\item $\mathcal{N}'(\bm{x}) \neq \emptyset$: in this case, we take each vertex from this set and add it to the input list, if does not belong to any of the orbits already explored, to maintain the list sorted by the number of zeros. Then, we start again from step~4.
\end{enumerate}
\item We continue iteratively until either $T_{tot}$ or $M$ are reached.
\end{enumerate}

Note that, thanks to Lemma \ref{lemma:neigh}, we do not need to pivot twice on isomorphic vertices, as they share the same neighborhood.
In the remainder of this section, we show step-by-step how the algorithm described in Section \ref{sec:four} works in practice from $n = 5$ to $ n = 6$.
First, we need a starting vertex. 
To do so, we recover it using the breaking loop procedure from $n = 5$.

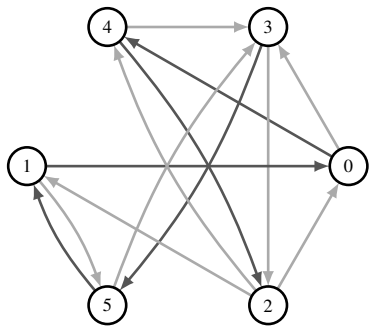
\begin{figure}[!t]
\centering
\begin{tikzpicture}[scale=0.75]
\Vertex[x=5.850,y=3.000,size=0.5,color=white,opacity=1,label=0,fontcolor=black]{0}
\Vertex[x=4.425,y=5.468,size=0.5,color=white,opacity=1,label=3,fontcolor=black]{3}
\Vertex[x=1.575,y=5.468,size=0.5,color=white,opacity=1,label=4,fontcolor=black]{4}
\Vertex[x=0.150,y=3.000,size=0.5,color=white,opacity=1,label=1,fontcolor=black]{1}
\Vertex[x=1.575,y=0.532,size=0.5,color=white,opacity=1,label=5,fontcolor=black]{5}
\Vertex[x=4.425,y=0.532,size=0.5,color=white,opacity=1,label=2,fontcolor=black]{2}
\Edge[,lw=1.0,color={170.000085,170.000085,170.000085},Direct,RGB](0)(3)
\Edge[,lw=1.0,color={84.99991499999999,84.99991499999999,84.99991499999999},Direct,RGB](0)(4)
\Edge[,lw=1.0,color={170.000085,170.000085,170.000085},Direct,RGB](3)(2)
\Edge[,lw=1.0,color={84.99991499999999,84.99991499999999,84.99991499999999},Direct,RGB, bend=10](3)(5)
\Edge[,lw=1.0,color={84.99991499999999,84.99991499999999,84.99991499999999},Direct,RGB, bend=10](4)(2)
\Edge[,lw=1.0,color={170.000085,170.000085,170.000085},Direct,RGB](4)(3)
\Edge[,lw=1.0,color={84.99991499999999,84.99991499999999,84.99991499999999},Direct,RGB](1)(0)
\Edge[,lw=1.0,color={170.000085,170.000085,170.000085},Direct,RGB, bend=10](1)(5)
\Edge[,lw=1.0,color={84.99991499999999,84.99991499999999,84.99991499999999},Direct,RGB, bend=10](5)(1)
\Edge[,lw=1.0,color={170.000085,170.000085,170.000085},Direct,RGB, bend=10](5)(3)
\Edge[,lw=1.0,color={170.000085,170.000085,170.000085},Direct,RGB](2)(0)
\Edge[,lw=1.0,color={170.000085,170.000085,170.000085},Direct,RGB](2)(1)
\Edge[,lw=1.0,color={170.000085,170.000085,170.000085},Direct,RGB, bend=10](2)(4)
\end{tikzpicture}%     without .tex extension
\caption{Starting point for the pivoting algorithm for $n = 6$.\label{fig:x0_6}
}
\end{figure}

This gives us a representative for 6 different orbits and, among them, we choose the one with the smallest number of zeros.
Figure \ref{fig:x0_6} reports its support graph.
Such vertex has entries equal to $\frac{k}{3}$, $k \in \{0, \ldots, 5\}$.
We start by pivoting from this vertex and collecting all its neighbors.
Then, among its neighbors, we move to the one having \rev{the smallest} number of zeros, and we proceed iteratively.
Different from the pivoting pipeline presented in \cite{elliott2008integrality}, we can recover at least one representative for each orbit.

%%%%%%%%%%%%%%%%%%%%%%%%%
\begin{algorithm}
\DontPrintSemicolon
\SetAlgoLined
%\SetAlgoNoLine
\KwInput{$n$, Number of nodes of the ATSP instance under study}
\KwInput{$R := \{\bm{x}_i\}_{i\in I}$, list of vertices available for $n-1$}
\KwInput{$M$, Maximum number of iteration\rev{s} of algorithm}
\KwInput{$T_{tot}$, time limit for the whole iteration\rev{s}}
\KwInput{$T_{it}$, timelimit for the single iteration}
\KwOutput{$R'$, a collection of non isomorphic vertices of $\pasep{n}$}
$R'$ = Extend($R$) \;
$R'$ = SortByNumberOfZeros($R'$) \;
$ i = 0$ \;
$\bm{x}_0 = R'[i]$ \;
$ Ps = [ ] $ \;
\While{$i < M$ and time $< T_{tot}$}{
\If {$\bm{x}_0$ is not isomorphic to any vertex in $Ps$}{
$Ps$.append($\bm{x}_0$) \;
$\mathcal{N}'(\bm{x_0}) =$ Pivoting($\bm{x}_0,T_{it}$) \;
\For{$\bm{y} \in \mathcal{N}(\bm{x_0})$}{
\If {$\bm{y}$ is not isomorphic to any vertex in $R'$}{
$R'$.insert($\bm{y}$)
}
} 
$i \gets i+1$\;
$Ps$.append($R'[i]$) \;
}
}  
\Return $R'$
\caption{Generating vertices of $\pasep{n}$ \rev{via} the explore-exploit algorithm.\label{alg:expl}}
\end{algorithm}
%%%%%%%%%%%%%%%%%%%%%%%%%%%%%

%%%%%%%%%%%%%%%%%%%%%%%%%%%%%%%%%%%%%%
% Results
%%%%%%%%%%%%%%%%%%%%%%%%%%%%%%%%%%%%%%%% 
\section{Computational results}
\label{sec:results}
This section shows the results of our experiments, focusing on the loop-breaking procedure's impact on expanding loops and the effect on the integrality gap. Additionally, it examines the neighborhood and the stabilizer of vertices for small values of $n$, focusing on the relation with the integrality gap.

\paragraph{Implementation details} All the experiments run on an \rev{$\times$}86\_64 architecture with a 13th Gen Intel(R) Core(TM) i5-13600 processor, offering 20 CPU cores with 32-bit and 64-bit operation modes. The pivoting algorithm is implemented in C++, while the iterative procedure is in Python. 
We use the Eigen library \cite{eigenweb} to deal with matrices and MPFR \cite{fousse2007mpfr} for the operations in multiple precision. %as we have observed that most of our matrices are ill-conditioned.
The integrality gap is computed by solving the linear program \eqref{gap:obj_func}--\eqref{const:gap_3} using Gurobi v9.5.0 \cite{gurobi}.

\subsection{Combining the symmetry breaking pivoting and the $\lambda$-loop breaking procedure}
\label{sec:impact_ll}
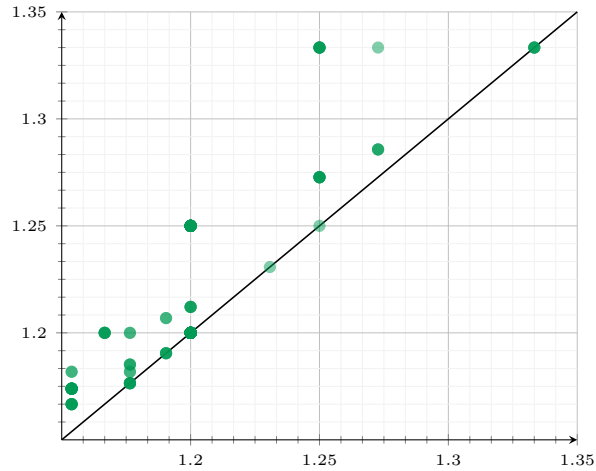
\begin{figure}[t!]
\centering
% This file was created with tikzplotlib v0.10.1.
\begin{tikzpicture}

\definecolor{mygreen}{HTML}{009B55}
\definecolor{steelblue31119180}{RGB}{31,119,180}

\begin{axis}[
    xmin=1.15,xmax=1.35,
    ymin=1.15,ymax=1.35,
    grid=both,
    grid style={line width=.1pt, draw=gray!10},
    major grid style={line width=.2pt,draw=gray!50},
    axis lines=middle,
    minor tick num=5,
    ticklabel style={font=\tiny,fill=white},
]
\addplot [semithick, mygreen, opacity=0.5, mark=*, mark size=2, mark options={solid}, only marks]
table {%
1.19047619047619 1.19047619047619
1.19047619047619 1.20689655172414
1.19047619047619 1.19047619047619
1.19047619047619 1.20689655172414
1.19047619047619 1.19047619047619
1.19047619047619 1.19047619047619
1.2 1.2
1.2 1.25
1.2 1.2
1.2 1.25
1.2 1.2
1.2 1.25
1.17647058823529 1.18181818181818
1.17647058823529 1.2
1.17647058823529 1.18518518518519
1.17647058823529 1.17647058823529
1.23076923076923 1.23076923076923
1.17647058823529 1.17647058823529
1.17647058823529 1.17647058823529
1.17647058823529 1.18518518518519
1.17647058823529 1.17647058823529
1.17647058823529 1.18518518518519
1.17647058823529 1.18181818181818
1.17647058823529 1.2
1.2 1.21212121212121
1.2 1.25
1.2 1.21212121212121
1.2 1.25
1.2 1.2
1.2 1.25
1.2 1.2
1.2 1.21212121212121
1.2 1.25
1.2 1.2
1.2 1.25
1.2 1.25
1.2 1.21212121212121
1.2 1.25
1.2 1.2
1.17647058823529 1.17647058823529
1.17647058823529 1.17647058823529
1.2 1.2
1.2 1.2
1.2 1.2
1.15384615384615 1.17391304347826
1.15384615384615 1.16666666666667
1.15384615384615 1.16666666666667
1.2 1.2
1.2 1.2
1.2 1.2
1.15384615384615 1.17391304347826
1.15384615384615 1.16666666666667
1.2 1.2
1.2 1.2
1.2 1.2
1.15384615384615 1.17391304347826
1.15384615384615 1.17391304347826
1.15384615384615 1.17391304347826
1.15384615384615 1.18181818181818
1.15384615384615 1.17391304347826
1.15384615384615 1.17391304347826
1.2 1.25
1.2 1.2
1.15384615384615 1.16666666666667
1.15384615384615 1.17391304347826
1.15384615384615 1.16666666666667
1.2 1.2
1.2 1.2
1.15384615384615 1.17391304347826
1.15384615384615 1.18181818181818
1.2 1.2
1.2 1.2
1.2 1.2
1.2 1.2
1.2 1.2
1.2 1.2
1.15384615384615 1.17391304347826
1.15384615384615 1.17391304347826
1.2 1.2
1.2 1.2
1.2 1.2
1.2 1.2
1.15384615384615 1.17391304347826
1.15384615384615 1.16666666666667
1.2 1.2
1.2 1.2
1.2 1.2
1.2 1.2
1.16666666666667 1.2
1.16666666666667 1.2
1.2 1.2
1.2 1.2
1.2 1.2
1.2 1.2
1.2 1.2
1.15384615384615 1.17391304347826
1.15384615384615 1.17391304347826
1.2 1.2
1.2 1.2
1.2 1.2
1.2 1.2
1.2 1.2
1.2 1.25
1.2 1.2
1.2 1.2
1.2 1.2
1.2 1.2
1.2 1.2
1.2 1.25
1.2 1.25
1.2 1.2
1.2 1.2
1.2 1.2
1.2 1.2
1.2 1.25
1.2 1.2
1.2 1.2
1.2 1.2
1.2 1.2
1.2 1.2
1.2 1.25
1.2 1.2
1.2 1.2
1.2 1.25
1.2 1.25
1.2 1.25
1.25 1.33333333333333
1.25 1.33333333333333
1.25 1.27272727272727
1.33333333333333 1.33333333333333
1.33333333333333 1.33333333333333
1.33333333333333 1.33333333333333
1.33333333333333 1.33333333333333
1.25 1.27272727272727
1.25 1.27272727272727
1.25 1.25
1.27272727272727 1.28571428571429
1.27272727272727 1.28571428571429
1.27272727272727 1.33333333333333
1.27272727272727 1.28571428571429
1.2 1.25
1.2 1.2
1.2 1.25
1.2 1.25
1.2 1.25
1.16666666666667 1.2
1.16666666666667 1.2
1.16666666666667 1.2
1.2 1.2
1.2 1.25
1.25 1.33333333333333
1.25 1.33333333333333
1.2 1.2
1.2 1.2
1.2 1.2
1.2 1.2
1.2 1.2
1.2 1.25
1.2 1.25
1.2 1.25
1.25 1.27272727272727
1.25 1.33333333333333
1.2 1.25
1.25 1.27272727272727
1.25 1.27272727272727
1.25 1.33333333333333
1.2 1.2
1.2 1.25
1.2 1.25
};
\addplot [semithick, mark=--, black]
table {%
1.15 1.15
1.35 1.35
};
\end{axis}

\end{tikzpicture}%     without .tex extension
\caption{On the $x$ axis, we represent the integrality gap of the vertices for $n = 6$ that have a $\lambda$-loop. For each of these vertices, we plot the integrality gap of the vertices obtained by the $\lambda$-loop breaking procedure on the $y$ axis. The fact that all points lie (non-strictly) above the line $y = x$ implies that the $\lambda$-loop breaking procedure is highly effective in increasing the integrality gap.}
\label{fig:gap_6_7}
\end{figure}

\begin{figure*}
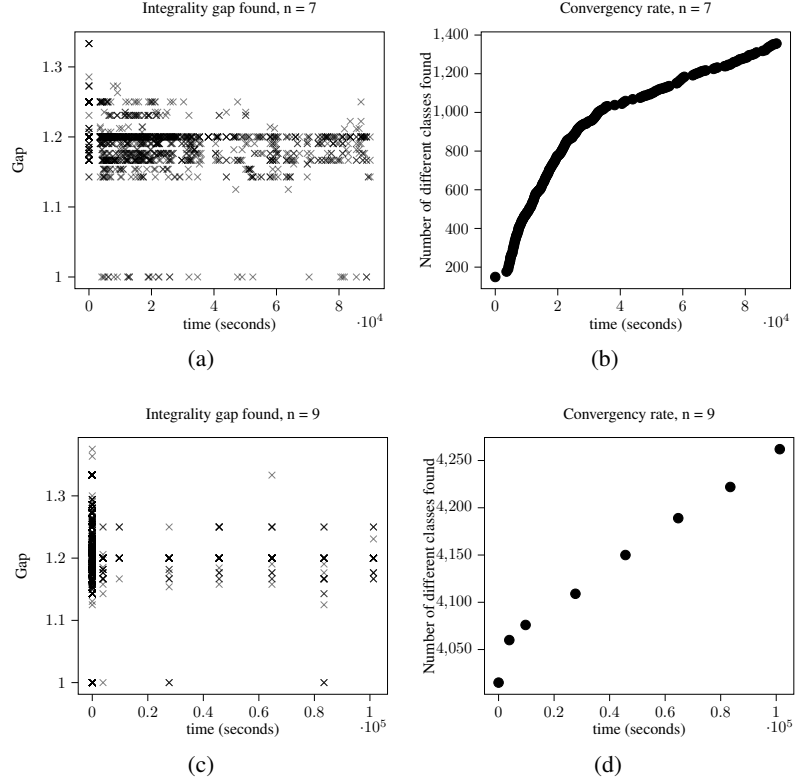

\centering
\subfloat[]{\scalebox{.6}{\input{include/gap_7}}
}~
\subfloat[]{
\scalebox{.6}{% This file was created with tikzplotlib v0.10.1.
\begin{tikzpicture}

\definecolor{darkgray176}{RGB}{176,176,176}

\begin{axis}[
tick align=outside,
tick pos=left,
title={Convergency rate, n = 7},
x grid style={darkgray176},
xlabel={time (seconds)},
xmin=-4498.50708400011, xmax=94468.6487640023,
xtick style={color=black},
y grid style={darkgray176},
ylabel={Number of different classes found},
ymin=88.65, ymax=1416.35,
ytick style={color=black}
]
\addplot [semithick, black, mark=*, mark size=3, mark options={solid}, only marks]
table {%
0 149
3623.62775802612 177
4110.3746778965 189
4260.9481408596 202
4420.36709070206 214
4601.02176046371 224
4751.40768504143 234
4812.36377668381 244
4890.28276705742 251
5044.32033348083 258
5204.67083144188 266
5362.42334151268 270
5520.54203009605 277
5670.97720217705 292
5821.06025004387 303
5973.41048765182 311
6125.73208332062 323
6279.92492508888 335
6432.49423933029 343
6590.76696753502 352
6749.79187464714 361
6823.53524303436 365
7009.5178706646 371
7168.43085026741 379
7327.47583627701 386
7411.97217559814 394
7480.46874451637 401
7640.40862679482 407
7800.85223484039 412
7957.43560528755 417
8140.55434370041 424
8297.61156368256 430
8457.10102677345 438
8627.94753479958 442
8785.29368019104 446
8970.11534404755 451
9129.0098965168 457
9287.77304935455 463
9452.6555249691 466
9607.20842051506 467
9749.98046684265 473
9905.0391125679 474
10079.7773628235 481
10252.6440000534 485
10419.7768573761 486
10583.0093786716 491
10747.238186121 495
10919.721716404 504
11083.8852031231 507
11248.9024963379 508
11422.16345191 514
11595.821228981 521
11772.731965065 530
11947.9677622318 535
12124.7045538425 542
12292.9801366329 547
12464.4178781509 558
12622.811157465 565
12795.4399380684 573
12955.3740849495 580
13124.0805869102 582
13292.1544966698 584
13498.3311822414 588
13671.5086283684 592
13852.0182025433 596
14025.1466715336 598
14200.8785293102 600
14381.4331715107 603
14583.3357005119 609
14777.575820446 615
14950.4361915588 624
15132.2183289528 631
15336.7384774685 636
15533.5721879005 641
15716.1889414787 646
15892.0968427658 655
16067.5988526344 661
16243.5217673779 667
16426.279638052 671
16609.2134225368 679
16802.7099084854 682
16984.8102099895 688
17168.3378949165 697
17263.5581576824 702
17373.4806442261 705
17560.6013095379 709
17750.302312851 714
17939.249874115 720
18128.9042658806 726
18328.0787448883 729
18515.3971097469 736
18697.413269043 740
18880.3999483585 744
19068.2943148613 748
19258.8801751137 757
19450.7837901115 761
19640.3424651623 770
19843.0347530842 774
20038.2859251499 777
20237.4372348785 779
20436.859172821 781
20632.4195249081 785
20827.6721737385 790
21030.7116594315 793
21227.8294355869 802
21430.1128423214 805
21625.9339444637 812
21819.9767522812 818
22013.9995789528 824
22217.7688624859 830
22417.7424612045 839
22623.6541411877 843
22831.7501022816 849
23039.6564948559 856
23460.5787675381 858
23681.3032233715 862
23901.9010589123 865
24118.6352925301 867
24564.5388314724 869
24774.7376270294 874
24982.1304719448 882
25402.237912178 888
25615.8930718899 892
25830.5405023098 897
26049.696190834 904
26268.4340598583 910
26486.2216169834 915
26708.3831276894 918
26928.0850820541 922
27147.5634527206 926
27369.7119276524 929
27593.7061259747 931
27814.0205667019 935
28261.7466037273 938
28482.8491053581 941
28936.7462289333 942
29162.5576252937 943
29719.3891558647 952
30274.1712024212 955
30833.9779167175 958
31065.8827176094 962
31293.4910581112 965
31530.2473340034 969
31763.8727362156 974
31992.8756537437 979
32229.5453150272 984
32462.9731259346 988
32710.3342840672 996
32942.7538094521 1000
33178.9232897758 1004
33723.7716302872 1010
34257.7198183537 1014
34773.4514975548 1019
35339.6332995892 1029
35872.5817973614 1032
38208.9744899273 1038
40528.5118124485 1045
41093.351474762 1053
41670.9525673389 1059
43971.2402610779 1068
46356.9674050808 1076
46968.8296451569 1080
47559.0558133125 1085
47813.7692708969 1086
48436.2077596188 1088
48996.9921247959 1093
49582.0992267132 1094
50148.219653368 1100
50692.6795363426 1103
51254.2118067741 1110
51841.1021864414 1114
52438.9677960873 1119
52973.5158042908 1122
53523.0227165222 1123
54052.8709681034 1127
54604.0346763134 1129
55137.7025206089 1133
55672.3208680153 1136
57962.4355552197 1151
58227.7708508968 1152
58791.0988311768 1161
59367.6194655895 1169
59932.7762935162 1178
60491.4831907749 1184
63204.9103312492 1192
63779.5155923367 1197
64355.0725479126 1201
64951.2777383327 1204
65510.5202765465 1208
66111.0527014732 1213
66715.6100251675 1215
67356.0343296528 1217
69841.9711995125 1225
70426.1139798164 1229
71033.8126168251 1232
73357.0695571899 1239
73924.3973705769 1240
74508.4111151695 1246
75121.2040476799 1247
75706.4550023079 1253
76319.6693344116 1260
76936.3518366814 1262
77542.5342755318 1267
78156.1086807251 1274
78765.1082391739 1276
79342.7483108044 1279
79961.822861433 1282
80566.4800081253 1285
81161.5541305542 1295
81741.511377573 1297
82347.0188312531 1299
82960.4340615273 1303
83590.975673914 1310
85565.7339966297 1318
86432.9614524841 1322
87042.5271704197 1332
87633.875957489 1339
88207.5965340137 1344
88811.2778708935 1351
89380.2864973545 1353
89970.1416800022 1356
};
\end{axis}

\end{tikzpicture}}
}

\subfloat[]{
\scalebox{.6}{\input{include/gap_9}}
}~
\subfloat[]{
\scalebox{.6}{% This file was created with tikzplotlib v0.10.1.
\begin{tikzpicture}

\definecolor{darkgray176}{RGB}{176,176,176}

\begin{axis}[
tick align=outside,
tick pos=left,
title={Convergency rate, n = 9},
x grid style={darkgray176},
xlabel={time (seconds)},
xmin=-5063.74556715489, xmax=106338.656910253,
xtick style={color=black},
y grid style={darkgray176},
ylabel={Number of different classes found},
ymin=4002.65, ymax=4274.35,
ytick style={color=black}
]
\addplot [semithick, black, mark=*, mark size=3, mark options={solid}, only marks]
table {%
0 4015
3869.92464876175 4060
9731.5111989975 4076
27683.5218601227 4109
45698.2415828705 4150
64714.7259676456 4189
83458.6374137402 4222
101274.911343098 4262
};
\end{axis}

\end{tikzpicture}}
}
\caption{Left: time vs integrality gap found. The darker the $\times$, the bigger the number of vertices found leading to that integrality gap. We note that the highest integrality gap is found at the very beginning of the procedure. Right: Time vs Number of classes of isomorphism.  We note that for $n = 7$ the algorithm continuously finds new and new vertices; for $n = 9$, new vertices are more and more rare.}\label{fig:gap6}
\end{figure*}

\rev{The major strength of our algorithm} is that\rev{,} compared to the Pivoting algorithm introduced by \cite{elliott2008integrality}, for $n = 6$, we manage to quickly recover at least one representative for all orbits.
Another strength is the impact of the $\lambda$-loop breaking procedure in quickly identifying vertices with a large integrality gap. Figure \ref{fig:gap_6_7} illustrates this phenomenon in the transition from $n = 6$ to $n = 7$. Experimentally, we observe that the integrality gap of the vertices obtained through $\lambda$-loop breaking is always greater than or equal to the starting one. 
Another strength is that the $\lambda$-loop breaking procedure alone is capable of finding the vertex with the maximum integrality gap for each $n$. For example, in Figure~\ref{fig:gap6}, it can be seen that the vertex with the maximum integrality gap was obtained at iteration 0, implying that the rest of the search, although leading to many points with large integrality gaps, is subordinated to what is obtained at step 0.

Among the weaknesses, we observe that although the solutions we found with the $\lambda$-loop loops procedure are associated with large integrality gaps, they have a lot of zeros. In fact, moving from $n$ to $n+1$, we switch from considering points from dimension $n(n-1)$ to points in $\erre^{(n+1)n}$, adding hence $2n$ entries. All of them are 0, but 2. So we add $2n - 2$ zeros to our vertices. Hence, it is hard to explore the full neighborhood, due to the great amount of feasible basis.
Lastly, we were able to push our pivoting procedure until $n = 11$. After that, it becomes infeasible to even partially enumerate at least one neighborhood of the vertex.
Figure \ref{fig:gap6} illustrates how the duration of each iteration increases and how the number of new orbits found decreases with each iteration.

%%%%%%%%%%%%%%%%%%

\subsection{Neighborhood exploration}
This section is devoted to studying the neighborhood of some vertices of $n \in \{4, 5, 6\}$ and deducing local properties.
Before diving into the details, let us recall the definition of \emph{polyhedral graph}.

\begin{definition}[Polyhedral graph]
	A polyhedral graph is an undirected graph formed from the vertices and edges of a convex polyhedron.
\end{definition}

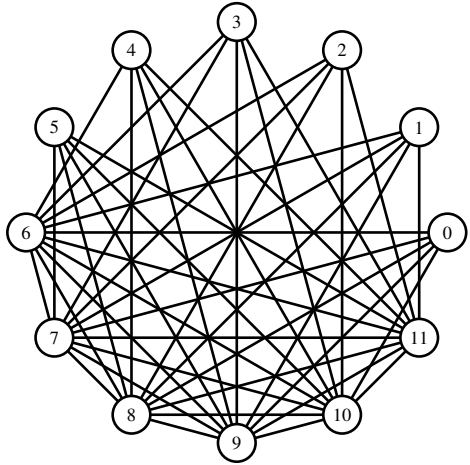
\begin{figure}[t!]
\centering
\begin{tikzpicture}
\Vertex[x=5.800,y=3.000,size=0.5,color=white,opacity=1,label=0,fontcolor=black]{0}
\Vertex[x=5.425,y=4.400,size=0.5,color=white,opacity=1,label=1,fontcolor=black]{1}
\Vertex[x=4.400,y=5.425,size=0.5,color=white,opacity=1,label=2,fontcolor=black]{2}
\Vertex[x=3.000,y=5.800,size=0.5,color=white,opacity=1,label=3,fontcolor=black]{3}
\Vertex[x=1.600,y=5.425,size=0.5,color=white,opacity=1,label=4,fontcolor=black]{4}
\Vertex[x=0.575,y=4.400,size=0.5,color=white,opacity=1,label=5,fontcolor=black]{5}
\Vertex[x=0.200,y=3.000,size=0.5,color=white,opacity=1,label=6,fontcolor=black]{6}
\Vertex[x=0.575,y=1.600,size=0.5,color=white,opacity=1,label=7,fontcolor=black]{7}
\Vertex[x=1.600,y=0.575,size=0.5,color=white,opacity=1,label=8,fontcolor=black]{8}
\Vertex[x=3.000,y=0.200,size=0.5,color=white,opacity=1,label=9,fontcolor=black]{9}
\Vertex[x=4.400,y=0.575,size=0.5,color=white,opacity=1,label=10,fontcolor=black]{10}
\Vertex[x=5.425,y=1.600,size=0.5,color=white,opacity=1,label=11,fontcolor=black]{11}
\Edge[,lw=1.0,color=black](0)(10)
\Edge[,lw=1.0,color=black](0)(6)
\Edge[,lw=1.0,color=black](0)(9)
\Edge[,lw=1.0,color=black](0)(8)
\Edge[,lw=1.0,color=black](0)(7)
\Edge[,lw=1.0,color=black](1)(11)
\Edge[,lw=1.0,color=black](1)(8)
\Edge[,lw=1.0,color=black](1)(6)
\Edge[,lw=1.0,color=black](1)(7)
\Edge[,lw=1.0,color=black](1)(9)
\Edge[,lw=1.0,color=black](2)(8)
\Edge[,lw=1.0,color=black](2)(7)
\Edge[,lw=1.0,color=black](2)(11)
\Edge[,lw=1.0,color=black](2)(6)
\Edge[,lw=1.0,color=black](2)(10)
\Edge[,lw=1.0,color=black](3)(9)
\Edge[,lw=1.0,color=black](3)(10)
\Edge[,lw=1.0,color=black](3)(7)
\Edge[,lw=1.0,color=black](3)(11)
\Edge[,lw=1.0,color=black](3)(6)
\Edge[,lw=1.0,color=black](4)(6)
\Edge[,lw=1.0,color=black](4)(9)
\Edge[,lw=1.0,color=black](4)(8)
\Edge[,lw=1.0,color=black](4)(10)
\Edge[,lw=1.0,color=black](4)(11)
\Edge[,lw=1.0,color=black](5)(7)
\Edge[,lw=1.0,color=black](5)(11)
\Edge[,lw=1.0,color=black](5)(10)
\Edge[,lw=1.0,color=black](5)(9)
\Edge[,lw=1.0,color=black](5)(8)
\Edge[,lw=1.0,color=black](6)(11)
\Edge[,lw=1.0,color=black](6)(7)
\Edge[,lw=1.0,color=black](6)(10)
\Edge[,lw=1.0,color=black](6)(8)
\Edge[,lw=1.0,color=black](6)(9)
\Edge[,lw=1.0,color=black](7)(9)
\Edge[,lw=1.0,color=black](7)(8)
\Edge[,lw=1.0,color=black](7)(10)
\Edge[,lw=1.0,color=black](7)(11)
\Edge[,lw=1.0,color=black](8)(9)
\Edge[,lw=1.0,color=black](8)(10)
\Edge[,lw=1.0,color=black](8)(11)
\Edge[,lw=1.0,color=black](9)(10)
\Edge[,lw=1.0,color=black](9)(11)
\Edge[,lw=1.0,color=black](10)(11)
\end{tikzpicture}%     without .tex extension
\caption{Polyhedral graph for $\pasep{4}$.}
\label{fig:adj_4}
\end{figure}

\rev{For \( n = 4 \), we used Polymake \cite{assarf2017computing} to enumerate all the vertices of \( \pasep{4} \). This task was previously performed in \cite{elliott2008integrality} using PORTA \cite{christof2009porta}.} More specifically, we have 12 vertices in total and two orbits. With this small number of vertices and orbits is hence easy to exhaustively study the neighborhood of each point.
Figure \ref{fig:adj_4} show the polyhedral graph of $\pasep{4}$. 
Nodes from 0 to 5 represent the non-integer vertices, while nodes from 6 to 11 represent the tours.
We can observe that each non-integer vertex has among its adjacent vertices always an integer one.
Interestingly, each tour is connected to all the vertices but one.

\begin{table}[t!]
    \caption{Orbit structure for $n = 4$. Columns: cardinality of the orbit, type of components, frequency of each component, and integrality gap attained at the elements of that orbit.}\label{tab:orbit_4}\centering
    \medskip
\begin{tabular}{lrrrrr}
\hline
$\vert O_{\bm{x}} \vert$ & \multicolumn{2}{c}{Components} & \multicolumn{2}{c}{Frequencies} & IG  \\
\hline
6      & 0        & \rev{1/2}       & 4           & 8         & $6/5$ \\
6      & 0        & 1       & 8           & 4         & 1\\
\hline
\end{tabular}
\end{table}

For $ n = 5$, we recover 384 vertices, as already done in \cite{elliott2008integrality}.
We choose a representative for each class and study its neighborhood. 
This can be done w.l.o.g. thanks to Lemma \ref{lemma:neigh}.
As expected, the number of neighbors is related to the degeneracy of the vertex, that is\rev{, the} number of zeros among both arc and slack variables.
Interestingly, each vertex has at least one adjacent representative for each equivalence class.
Furthermore, the vertices (d) and (e) also have a representative of themselves among the neighbors.

\begin{table}[t!]
    \caption{Orbit structure for $n = 5$. Columns: Label the orbits according to Figure \ref{fig:n5repr}, Cardinality of the orbit, type of components, frequency of each component,  integrality gap attained at the elements of that orbit, number of tight sets, neighborhood size, and stabilizer.}\label{tab:orbit_5}
\centering\medskip
\begin{tabular}{lrrrrrrrrrrr}
\hline\noalign{\smallskip}
&$\vert O_{\bm{x}} \vert $ & \multicolumn{3}{c}{Components} & \multicolumn{3}{c}{Frequencies} & IG   & tight sets & $\vert \mathcal{N}(\bm{x}) \vert $ & $G_{\bm{x}}$\\
\noalign{\smallskip}\hline\noalign{\smallskip}
(a) & 60     & 0          & $6/5$        &            & 10          & 10 &           & $5/4$  & 6 & 28 & $\langle (0 \; 1)(4 \; 3)\rangle$\\
(b) & 120    & 0          & $\frac{1}{3}$ & $\frac{2}{3}$ & 9            & 7  & 4         & $6/5$ & 4 & 20 & $\langle id \rangle$ \\
(c) & 60     & 0          & $\frac{1}{2}$        &            & 10          & 10 &           & $6/5$ & 4 & 28 & $\langle (0 \; 4)(1 \; 3)\rangle$\\ 
(d) & 120    & 0          & 1          & $\frac{1}{2}$        & 11          & 1  & 8         & $6/5$ & 4 & 23 & $\langle id \rangle$ \\
(e) & 24     & 0          & 1          &            & 15          & 5  &           & 1  & 10 & 148 & $\langle (0 \; 1 \; 2\; 3 \; 4) \rangle$ \\
\noalign{\smallskip}\hline
\end{tabular}
\end{table}

In the case of $n=6$, we have 90 orbits. The orbit structure can be found in Table \ref{tab:orbit_6_1} and \rev{Table} \ref{tab:orbit_6_2}, in the appendix.
Neighborhoods of these vertices cannot be exhaustively explored. When the number of zeros is close to 18, the exhaustive listing of the whole neighborhood quickly leads to an out-of-memory error.

Note that, from our preliminary test, we state the following conjecture:
\begin{conj}
For each $\bm{x}$ vertex, $\mathcal{N}(\bm{x})$ contains at least one tour.
\end{conj}
If this conjecture were true, solely pivoting on \rev{an} integer vertex would lead to the full V-description of the polytope:
By using Lemma \ref{lemma:neigh}, the neighborhood of an integer vertex will contain at least one representative of each orbit.
We can then list all the vertices by fully describing each orbit.
This can be done just \emph{theoretically}, as in practice, listing the full neighborhood of a vertex with just 6 nodes is \rev{computationally} infeasible.

\subsection{Symmetries and conjectured relation with the integrality gap}
\label{sec:results_symmetries}
Table \ref{tab:orbit_4}, \ref{tab:orbit_5}, \ref{tab:orbit_6_1}, \ref{tab:orbit_6_2} show the orbit structure of $n = 4, 5, 6$. In all cases, we observe that the maximum integrality gap has been attained at a one-half integer vertex with a relatively small orbit.

Table \ref{tab:orbit_4} reports the structure of the orbit we found for $ n = 4$. As there is only one half-integer solution, the maximum integrality gap is attained at that vertex. 

For $n = 5$, we have represented all the vertices in Figure \ref{fig:n5repr}.
The one attaining the highest integrality gap is vertex (a). 
Looking at Table \ref{tab:orbit_5}, we understand that \rev{vertices} (a) and (c) have the same number of non-zero and zero entries.
The main differences \rev{are} in the number of $\lambda$-loops and tight sets (see Definition \ref{def:ts}).
Tight sets are associated with slack variables $\overline{x}_S = 0$.
In particular, vertex (a) has 6 tight sets, while vertex (c) has only 4. 
For $n = 6$, there are 90 orbits, see Table\rev{s} \ref{tab:orbit_6_1} and \ref{tab:orbit_6_2}. 
The vertex attaining the maximum integrality gap is again half-integer, which in principle does not seem very different from other half-integer vertices having low integrality gap (See, e.g, the middle of Table \ref{tab:orbit_6_2}).
We observe that all the half-integer vertices have 10 tight sets. The two vertices \rev{having} two highest integrality gaps have instead a higher number of tight sets (Line\rev{s} 1 and 3 of Table \ref{tab:orbit_6_1}). More specifically, the two half-integer vertices maximizing the integrality gap have, respectively, 10 and 12 tight sets.

For these small number\rev{s} of nodes, we also explicitly compute the stabilizer. 
For $n = 4$, the stabilizer of any representative of the non-integer orbit is given by the vertices of the subcycles in the support graph. For instance, referring to Figure \ref{fig:esempio_4}, left, it holds $G_x=\langle(0 \; 3 \; 1 \; 2)\rangle$. 

All the stabilizers are trivial for $n = 5$, but the ones of the two half-integer vertices. Table \ref{tab:orbit_5} reports all the stabilizers explicitly computed. The vertices attaining the maximum gap are the ones whose stabilizer is isomorphic to  $\mathbb{Z}_2$ and, more specifically, it acts by swapping the extremes of the two cycles. For \rev{vertex}(c), the situation is analogous, by considering the two chained $\lambda$-loops $\overleftrightarrow{03}$ and $\overleftrightarrow{04}$ as one.

\begin{figure}[!t]
\centering
\begin{tikzpicture}
\Vertex[x=3.000,y=5.468,size=0.5,color=white,opacity=1,label=0,fontcolor=black]{0}
\Vertex[x=4.425,y=5.468,size=0.5,color=white,opacity=1,label=4,fontcolor=black]{4}
\Vertex[x=1.575,y=5.468,size=0.5,color=white,opacity=1,label=5,fontcolor=black]{5}
\Vertex[x=1.575,y=3.000,size=0.5,color=white,opacity=1,label=1,fontcolor=black]{1}
\Vertex[x=3.000,y=3.000,size=0.5,color=white,opacity=1,label=3,fontcolor=black]{3}
\Vertex[x=4.425,y=3.000,size=0.5,color=white,opacity=1,label=2,fontcolor=black]{2}
\Edge[,lw=1.0,color={127.5,127.5,127.5},Direct,RGB, bend=10](0)(4)
\Edge[,lw=1.0,color={127.5,127.5,127.5},Direct,RGB, bend=10](0)(5)
\Edge[,lw=1.0,color={127.5,127.5,127.5},Direct,RGB, bend=10](4)(0)
\Edge[,lw=1.0,color={127.5,127.5,127.5},Direct,RGB](4)(1)
\Edge[,lw=1.0,color={127.5,127.5,127.5},Direct,RGB, bend=10](5)(0)
\Edge[,lw=1.0,color={127.5,127.5,127.5},Direct,RGB](5)(2)
\Edge[,lw=1.0,color={127.5,127.5,127.5},Direct,RGB, bend=10](1)(3)
\Edge[,lw=1.0,color={127.5,127.5,127.5},Direct,RGB](1)(5)
\Edge[,lw=1.0,color={127.5,127.5,127.5},Direct,RGB, bend=10](3)(1)
\Edge[,lw=1.0,color={127.5,127.5,127.5},Direct,RGB, bend=10](3)(2)
\Edge[,lw=1.0,color={127.5,127.5,127.5},Direct,RGB, bend=10](2)(3)
\Edge[,lw=1.0,color={127.5,127.5,127.5},Direct,RGB](2)(4)
\end{tikzpicture}
 \caption{Vertex obtaining the maximum integrality gap for $ n = 6$.\label{fig:vertex_max_gap}}
\end{figure}
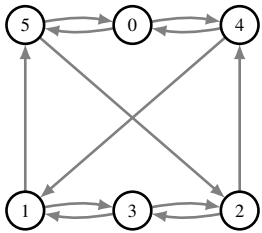

For $n = 6$, the vertex attaining the maximum integrality gap has been represented in Figure \ref{fig:vertex_max_gap}. The stabilizer of this vertex is isomorphic to $\mathbb{Z}_4$, more specifically \rev{it} is the group generated from $(0 3)(1 4 2 5).$ 
Even in this case, the stabilizer ``swaps'' the extreme nodes of the two chained $\lambda$-loops and the middle nodes 0 and 3.

For $n = 7$, our heuristic procedure does not recover all the vertices, hence we can do considerations only among the ones we were able to obtain, which are 1356 out of 3748. According to \cite{elliott2008integrality}, 5 different isomorphism classes attain the maximum gap of $\frac{4}{3}$: our heuristic finds all of them. As already observed in \cite{elliott2008integrality}, 3 out of 5 vertices have entries in $\{0, 0.5\}$ (from now on, we will denote vertices having entries in $\{0, 0.5\}$ as \emph{half-integer}) while the others in $\{0, 0.5, 1\}$ (\emph{integer-half-integer}). Even in this case, the number of tight sets in vertices maximizing the gap is high, namely, 12 and 16, although not maximal, as there exist pure half-integer vertices having 14 and 16 tight sets. A similar argument holds for half-integer vertices.

For $n = 8$, we identified 41 vertices with a maximum integrality gap of $\frac{4}{3}$. In \cite{elliott2008integrality}, there were 43 of such vertices, but one of them was not half-integer. Observe that their method was specifically tailored to find \emph{all} the half-integer vertices and \emph{potentially} other\rev{s}, whereas our approach allows for more flexibility. Specifically, we discovered 17 pure half-integer orbits, 16 half-\rev{integer} orbits, and 8 orbits with components in $\{0, 0.25, 0.5, 0.75\}$. While this suggests the possibility that non-half-integer vertices may also result in the maximum integrality gap, it appears to be an isolated case.

More specifically, for $n \geq 9$, the vertex attaining the maximum gap is always unique and a pure half-integer vertex. 
Based on our recent discussion, we have gathered the following empirical evidence:
\begin{itemize}
    \item Among the vertices achieving the maximum integrality gap, there is at least one half-integer \rev{vertex}.
    \item Furthermore, vertices with a large stabilizer and a large number of tight sets appear to maximize the integrality.
\end{itemize}

\subsection{New lower bound on the integrality gap}
\label{sec:new_lb}
Table \ref{tab:final_collaps_extend} presents the results of the combined symmetry-breaking pivoting and $\lambda$-loop breaking algorithm.
Starting from a vertex of $\pasep{6}$, 
our combined algorithm alternates the exploration of vertices of $\pasep{n}$ with the generation of a new vertex of $\pasep{n+1}$

With our combined algorithm, we can recompute all the lower bounds of $\alpha_n^{LB}$ up to $n=15$, and we compute newer lower bounds for $n \in \{16, 18, 20, 22\}$
Note that for $n$ odd, \cite{elliott2008integrality} introduces a family of ATSP instances having $\alpha_n=\frac{3k+1}{2k+1}$ where $n = 3 + 2(k+1)$, which gives new lower bounds for $n=17,19,21$.
However, all the lower bounds in \cite{elliott2008integrality} are obtained only by exploring half-integer vertices, while our procedure can generate non-half-integer vertices, thanks to the $\lambda$-loop breaking procedure. 
However, except for the case $n=8$, where we found a non-half-integer vertex maximizing the gap, such a maximum is always achieved in correspondence with a half-integer vertex.

To obtain Table \ref{tab:final_collaps_extend},  we proceed as follows: starting from a vertex explicitly given by \cite{elliott2008integrality} for $n = 18$, having an integrality gap of 1.5, we make each $\lambda$-loop \emph{collapse} to one single node, and check time by time if the so-obtained \emph{fesible} point is a vertex. 
Afterward, we expand each $\lambda$-loop, obtaining vertices for $ 19 \leq n \leq 22$.
These two procedures built Table \ref{tab:final_collaps_extend}, \rev{where the best known lower bounds for the integrality gap have been reported, as well as the new ones we got with our procedure. We can observe that} the lower bounds have been improved with respect to the state of the art for $n \in \{16, 17, 19, 20, 21, 22\}$.
We recall that, although our exploration allows different types of fractional vertices, the maximum gap is always attained on a half-integer vertex.

\begin{table}[t!]
\setlength{\tabcolsep}{5pt}
    \centering
    \caption{State of the art on the lower bounds $\alpha_n^{LB} \leq \alpha_n$ for ATSP, with $11 \leq n \leq 22$. In bold, the new best lower bounds are obtained with our approach. \rev{Column ``$n$'' reports the problem's dimension; Column ``From \cite{elliott2008integrality}'' provides the best-known lower bound for a given $n$ as established in the referenced paper. Column ``New'' presents the lower bounds we obtained.}}
    \label{tab:final_collaps_extend}
    \smallskip
\begin{tabular}{rr@{ }lr@{ }lcr@{ }lr@{ }l}
    \hline\noalign{\smallskip}
    $n$ & 
    \multicolumn{2}{c}{From \cite{elliott2008integrality}} &  
    \multicolumn{2}{c}{New} & 
    $n$ & 
    \multicolumn{2}{c}{From \cite{elliott2008integrality}} &  
    \multicolumn{2}{c}{New} \\
    \noalign{\smallskip}\hline\noalign{\smallskip}                     
    11 & $10/7$ & $(1.429)$ & $10/7$  & $(1.429)$  & 17 & $19/13$& $(1.461)$ & $\bm{55/37}$ &  $\bm{(1.486)}$ \\
    12 & $56/39$ & $(1.436)$ & $56/39$ & $(1.436)$ & 18 & $3/2$ & $(1.500)$& $3/2$ & $(1.500)$ \\
    13 & $13/9$ & $(1.444)$ & $13/9$ &  $(1.444)$ & 19 & $22/15$& $(1.466)$ & $\bm{3/2}$ & $\bm{(1.500)}$\\ 
    14 & $100/69$ & $(1.449)$ & $100/69$ &  $(1.449)$ & 20 & -&& $\bm{3/2}$ & $\bm{(1.500)}$ \\
    15 & $16/11$ & (1.454) & $16/11$ &  (1.454) & 21 & $25/17$& $(1.470)$ & $\bm{3/2}$ & $\bm{(1.500)}$ \\
    16 & -&& $\bm{28/19}$ & $\bm{(1.474)}$         & 22 & -&& $\bm{3/2}$ & $\bm{(1.500)}$ \\
    \noalign{\smallskip}\hline
    \end{tabular}
\end{table}

%%%%%%%%%%%%%%%%%%%%%%%%%%%%%%%
\subsection{Hard-ATSPLIB instances}\label{sec:hard_ATSPLIB}

Whenever we solve the Gap$(\bm x)$ problem to compute the maximum integrality gap for a given vertex, we \rev{generate an ATSP instance that could be challenging for ATSP solvers in practice}.
Hence, we \rev{propose 13 small hard ATSP instances}, which we share online\footnote{\url{https://github.com/eleonoravercesi/HardATSPLIB}}.
Several studies have examined the empirical hardness of the STSP concerning the integrality gap, such as \rev{\cite{hougardy2020hard,vercesi2023generation,zhong2021lower}}.
However, to the best of our knowledge, no such studies have been conducted on the ATSP.
Note that every time we compute the maximum possible integrality gap attained at a given vertex $\bm{x}$ by solving problem \eqref{gap:obj_func}--\eqref{const:gap_3}, we obtain a cost vector $\bm{c}$ associated to an ATSP instance having $\text{ATSP}(\bm{c}) = 1$.
Hence, we can evaluate the computational complexity of each instance as generated in this way.

A core question is how to evaluate complexity from a computational perspective.
\rev{Previous work by \cite{fischetti1997polyhedral} has shown promising results using a branch-and-cut algorithm that exploits facet-defining inequalities for ATSP. 
However, more recent studies \cite{fischetti2007exact,roberti2012models} suggest that Concorde, the \rev{state-of-the-art} solver for the STSP, can also be a competitive method for solving ATSP instances after suitably converting them into STSP. More specifically, it} is important to note that Concorde can only handle symmetric nonnegative and integer costs, but it is possible to transform any of the ATSP instances we obtained into an integer and non-negative STSP starting from the method proposed in \cite{jonker1983transforming}. 
First of all, we have observed that all the solutions we found solving \rev{$\text{Gap}(\bm{x})$} are rational, and hence it is possible to make the costs integer by multiplying all the entries by the common denominator.
Hence, without loss of generality, we can consider all the solutions of \rev{$\text{Gap}(\bm{x})$} as integer vectors.
Let $\overline{\bm{C}}=(\overline{c}_{ij})$ be a matrix derived from the cost vector as suggested in \rev{the paper of Jonker and Volgenant} \cite{jonker1983transforming}, namely:
\begin{equation}
\overline{c}_{ij} = \begin{cases}
c_{ij} \qquad i \neq j \\
-M \qquad i = j,
\end{cases}
\label{eq:transfATSPtoSTSP}
\end{equation}
\noindent where $M$ denotes a large positive number. Consider the matrix $\bm{U}$, where all entries are set to infinity. We construct the following $\mathbb{R}^{2n \times 2n}$ matrix:
\[ \tilde{\bm{C}}=\left[\begin{array}{ll}\bar{\bm{C}} & \bm{U} \\ \bm{U}& \bar{\bm{C}}\end{array}\right]. \]
Note that $ \tilde{\bm{C}}$ may contain negative costs, which we do not want, as Concorde only performs with \revv{non-negative} costs. Therefore, we shift all costs forward by $M$, namely making the minimum cost equal to 0. 
Unfortunately, in this framework, we only have \emph{premetrics}, as we lose the triangle inequality in every triple involving two of the original nodes and one ``doubled'' node.
However, since we have only performed an affine translation on each cost, the optimal tour does not change. The relationship between the original values of ATSP and STSP is hence
\begin{equation}
\atsp{c} = \stsp{\tilde{\bm{c}}} - nM.
\end{equation}

Hence, we transform each ATSP \rev{instance} into an STSP \rev{instance} as discussed above, solve each instance using Concorde, and record the computational time and the integrality gap. 

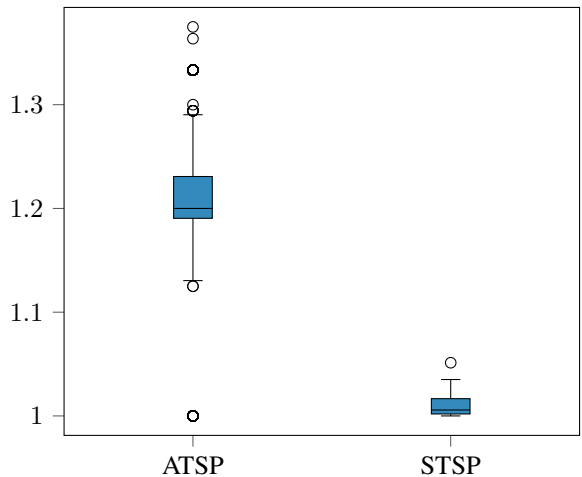
\begin{figure}[!t]
\centering
% This file was created with tikzplotlib v0.10.1.
\begin{tikzpicture}

\definecolor{dimgray85}{RGB}{85,85,85}
\definecolor{gainsboro229}{RGB}{229,229,229}
\definecolor{steelblue52138189}{RGB}{52,138,189}

\begin{axis}[
axis line style={black},
tick align=outside,
tick pos=left,
x grid style={white},
xmajorgrids,
xmin=0.5, xmax=2.5,
xtick style={color=dimgray85},
xtick={1,2},
xticklabels={ATSP,STSP},
y grid style={white},
ymajorgrids,
ymin=0.98125, ymax=1.39375,
ytick style={color=dimgray85}
]
\path [draw=black, fill=steelblue52138189]
(axis cs:0.925,1.19047619047619)
--(axis cs:1.075,1.19047619047619)
--(axis cs:1.075,1.23076923076923)
--(axis cs:0.925,1.23076923076923)
--(axis cs:0.925,1.19047619047619)
--cycle;
\addplot [black]
table {%
1 1.19047619047619
1 1.1304347826087
};
\addplot [black]
table {%
1 1.23076923076923
1 1.29032258064516
};
\addplot [black]
table {%
0.9625 1.1304347826087
1.0375 1.1304347826087
};
\addplot [black]
table {%
0.9625 1.29032258064516
1.0375 1.29032258064516
};
\addplot [black, mark=o, mark size=2, mark options={solid,fill opacity=0}, only marks]
table {%
1 1.00000022222227
1 1
1 1
1 1
1 1.12500014062502
1 1
1 1
1 1
1 1
1 1
1 1
1 1
1 1.00000066666711
1 1.000001000001
1 1.00000075000056
1 1.00000033333344
1 1.00000033333344
1 1.00000066666711
1 1.00000044444464
1 1.00000060000036
1 1.00000075000056
1 1.00000066666711
1 1
1 1
1 1.00000066666711
1 1
1 1.00000070967792
1 1.0000005238098
1 1.00000075000056
1 1.00000028571437
1 1.00000054545484
1 1.000001000001
1 1
1 1.00000050000025
1 1.00000050000025
1 1
1 1.00000040000016
1 1.00000072727326
1 1.00000083333403
1 1.00000076923136
1 1.00000072727326
1 1.00000080000064
1 1.00000081818249
1 1
1 1
1 1.00000066666711
1 1.00000033333344
1 1
1 1.125
1 1
1 1
1 1.33333333333333
1 1.33333333333333
1 1.33333333333333
1 1.33333333333333
1 1.33333333333333
1 1.33333333333333
1 1.33333333333333
1 1.33333333333333
1 1.33333333333333
1 1.33333333333333
1 1.33333333333333
1 1.33333333333333
1 1.33333333333333
1 1.33333333333333
1 1.33333333333333
1 1.33333333333333
1 1.33333333333333
1 1.33333333333333
1 1.33333333333333
1 1.33333333333333
1 1.33333333333333
1 1.33333333333333
1 1.33333333333333
1 1.33333333333333
1 1.33333333333333
1 1.33333333333333
1 1.33333333333333
1 1.33333333333333
1 1.33333333333333
1 1.33333333333333
1 1.33333333333333
1 1.33333333333333
1 1.33333333333333
1 1.33333333333333
1 1.33333333333333
1 1.33333333333333
1 1.33333333333333
1 1.33333333333333
1 1.33333333333333
1 1.33333333333333
1 1.33333333333333
1 1.33333333333333
1 1.33333333333333
1 1.33333333333333
1 1.33333333333333
1 1.33333333333333
1 1.33333333333333
1 1.33333333333333
1 1.33333333333333
1 1.33333333333333
1 1.33333333333333
1 1.33333333333333
1 1.33333333333333
1 1.33333333333333
1 1.33333333333333
1 1.33333333333333
1 1.33333333333333
1 1.33333333333333
1 1.33333333333333
1 1.33333333333333
1 1.33333333333333
1 1.33333333333333
1 1.33333333333333
1 1.33333333333333
1 1.33333333333333
1 1.33333333333333
1 1.33333333333333
1 1.33333333333333
1 1.33333333333333
1 1.33333333333333
1 1.33333333333333
1 1.33333333333333
1 1.33333333333333
1 1.33333333333333
1 1.33333333333333
1 1.33333333333333
1 1.33333333333333
1 1.33333333333333
1 1.33333333333333
1 1.33333333333333
1 1.33333333333333
1 1.33333333333333
1 1.33333333333333
1 1.33333333333333
1 1.375
1 1.36363636363636
1 1.33333333333333
1 1.33333333333333
1 1.33333333333333
1 1.33333333333333
1 1.33333333333333
1 1.33333333333333
1 1.33333333333333
1 1.33333333333333
1 1.33333333333333
1 1.33333333333333
1 1.33333333333333
1 1.33333333333333
1 1.33333333333333
1 1.33333333333333
1 1.33333333333333
1 1.33333333333333
1 1.33333333333333
1 1.33333333333333
1 1.33333333333333
1 1.33333333333333
1 1.33333333333333
1 1.29411764705882
1 1.29411764705882
1 1.29411764705882
1 1.29411764705882
1 1.29411764705882
1 1.29411764705882
1 1.33333333333333
1 1.3
1 1.33333333333333
};
\path [draw=black, fill=steelblue52138189]
(axis cs:1.925,1.00189109475357)
--(axis cs:2.075,1.00189109475357)
--(axis cs:2.075,1.01663223350523)
--(axis cs:1.925,1.01663223350523)
--(axis cs:1.925,1.00189109475357)
--cycle;
\addplot [black]
table {%
2 1.00189109475357
2 1
};
\addplot [black]
table {%
2 1.01663223350523
2 1.03508771929825
};
\addplot [black]
table {%
1.9625 1
2.0375 1
};
\addplot [black]
table {%
1.9625 1.03508771929825
2.0375 1.03508771929825
};
\addplot [black, mark=o, mark size=2, mark options={solid,fill opacity=0}, only marks]
table {%
2 1.05128205128205
};
\addplot [black]
table {%
0.925 1.20000024000005
1.075 1.20000024000005
};
\addplot [black]
table {%
1.925 1.00571020990422
2.075 1.00571020990422
};
\end{axis}

\end{tikzpicture}
  \caption{Distribution of the integrality gap from the ATSP to the STSP via the application of the Jonker-Volgenant-based procedure for $ n = 9$.}\label{fig:boxplot6}
\end{figure}

Regarding the integrality gap, we observe that we only have a relation between the value of $\atsp{\bm{c}}$ and $\stsp{\tilde{\bm{c}}}$, but nothing can be said \rev{\emph{a priori}} for $\ssep{\tilde{\bm{c}}}$ and $\asep{\bm{c}}$.
Remarkably, we observed that after the above-discussed procedure, the resulting instances exhibit a reduced integrality gap.
Figure \ref{fig:boxplot6} reports this information for $ n = 9$.

In terms of computational complexity, we are hence able to retrieve some hard-to-solve instances.
Table \ref{tab:comp_hard} reports the results for \rev{the hardest instance we get for every } $7 \leq n \leq 19$. We do \rev{six} different types of computation:
\begin{enumerate}
    \item We run Concorde as it is, and we record the runtime and the number of Branch \& Bound (B\&B) nodes. From now on, this would be called the ``standard'' setting, \rev{and in the Table is denoted with \texttt{-}}
    \item We add the flag \texttt{-d} that uses \rev{Depth-First} (DF) branching instead of Breadth-First (BF). Adding this flag can prevent Concorde from writing search nodes to files, which could not be a good idea for small instances.
    \item We add the flag \rev{\texttt{-C0}} to disable local cuts. For small \rev{STSP instances}, the computational overhead associated with generating and applying local cuts may outweigh the benefits they provide.
    \item We combine the flag\revv{s} \texttt{-d -C0} together.
    \item \rev{We run the algorithm presented in \cite{fischetti1992additive}, here denoted with \texttt{ATSP92}, where the authors present new lower bounds for the ATSP and implement an efficient B\&B exploiting them.} 
    \item \rev{We run the algorithm presented in \cite{fischetti2007exact}, here denoted with \texttt{FLT07}, where a Branch-and-Cut has been enhanced with the Fractionality Persistency criterion to improve the quality of the solution by prioritizing persistently fractional variables for branching.}
\end{enumerate}

\rev{Concorde is written in C, compiled with \texttt{gcc 13.3.0} and linked with CPLEX 12.9.0 \cite{cplex}. \texttt{ATSP92} and \texttt{FLT07} are written in Fortran and compiled with \texttt{gfortran-9}. The LP in \texttt{FLT07} is solved using CPLEX 12.9.0.} 

For comparison, we run the same set of experiments on the instances of the ATSPLIB \cite{reinelt1991tsplib} having fewer than 100 nodes. Results can be found in Table \ref{tab:comp_hard} and Table \ref{tab:comp_hard_ATSPLIB}.

First, we observe that the instances of the ATSPLIB are much easier for Concorde with respect to the ones of the HardATSPLIB. The instance $\texttt{ftv38}$ that is twice as large as our instance with $19$ nodes requires an order of magnitude less time. 
Overall, as expected, using DF and disabling local cuts greatly benefits the impact of Concorde. In the hardest instances, this leads to a faster runtime compared to \texttt{FLT07}.
Note that, overall, Concorde is on average slower than both \texttt{ATSP95} and \texttt{FLT07}, regardless of the setting. 
For the hard instances, although the optimal solution is often found at the root node via primal heuristics, a large number of branching nodes are generated. Each of these nodes yields only a slight improvement in the best lower bound, which results in long runtimes to prove optimality and a lot of branching nodes. For the instances of the ATSPLIB, the lower bound is of high quality, and as such, closing the branching tree is relatively fast. 

The dedicated methods \texttt{ATSP95} and \texttt{FLT07} appear more promising.
\texttt{ATSP95} relies on additive bounding techniques. This makes the computation per node extremely fast. However, the absence of stronger cuts causes a rapid growth in the number of branching nodes as the instance size increases. As a result, the number of branching nodes becomes large even for relatively small instances. For instance, in \texttt{p43}, the queue size exceeds the maximum limit of 1'500'000, which prevents the solver from returning the optimal solution. However, all the instances from the HardATSPLIB are successfully solved, although the computational effort is clearly higher: for $n = 17$, the solver requires 218'251 nodes and over one second, while \texttt{br17} is solved with only 22 nodes.

\texttt{FLT07} implements a dedicated branch-and-cut algorithm specifically designed for the ATSP. It outperforms Concorde on small instances and, although slightly slower on some larger ones \cite{fischetti2007exact}, it performs best overall on the ATSPLIB — both in terms of runtime and number of nodes. However, some interesting behavior appears when we move to the HardATSPLIB: for $n = 17, 18$, an error occurs during the root node solution phase, which prevents us from gathering complete performance data. For all other instances, the solver can prove optimality within a few seconds. Notably, the hard instance with $n = 16$ requires more than 8 seconds — nearly two orders of magnitude more than \texttt{kro124p} ($n = 100$) — and about 500 times more branching nodes.
%(``\texttt{PR_BUILD: ROOT ERROR}'').

\begin{sidewaystable}
%\sidewaystable
\begin{center}
\begin{minipage}{\textheight}
\footnotesize
\centering
\caption{Results of the runs of different solvers on the hard instances. The first column represents the number of nodes in the ATSP instances; \texttt{(-)} indicates that Concorde was run with the default settings; \texttt{-d} refers to running Concorde with the DF option enabled; \texttt{-C0} indicates that local cuts were disabled; \texttt{-d -C0} corresponds to running Concorde with both DF enabled and local cuts disabled. \texttt{ATSP92} reports the performance from \revv{\cite{fischetti1992additive}}, and \texttt{FLT07} shows the performance of the algorithm presented in \revv{\cite{fischetti2007exact}}. The last column shows the integrality gap (\texttt{IG}) for both instances—the original ATSP and the transformed STSP, obtained through the procedure described in Section \ref{sec:hard_ATSPLIB}. Gaps are rounded to the third decimal place for readability. \revv{Character ``-'' means that \texttt{FLT07} was not able to solve the root node. ``$^*$'' means that the queue size limit has been reached.}}\label{tab:comp_hard}
\begin{adjustbox}{scale=0.8}
\begin{tabular}{lrrrrrrrrrrrrrr}
\hline\noalign{\smallskip}
$n$ & \multicolumn{2}{c}{ \rev{-} } & \multicolumn{2}{c}{\rev{\texttt{-d}}} & \multicolumn{2}{c}{\rev{\texttt{-C0}}} & \multicolumn{2}{c}{\rev{\texttt{-C0 -d}}}  & \multicolumn{2}{c}{\rev{\texttt{ATSP92}}} & \multicolumn{2}{c}{\rev{\texttt{FLT07}}} & \multicolumn{2}{c}{\rev{\texttt{IG}}}\\
& \shortstack{Time\\(s)} & \shortstack{B\&B\\nds} & \shortstack{Time\\(s)} & \shortstack{B\&B\\nds} & \shortstack{Time\\(s)} & \shortstack{B\&B\\nds} & \shortstack{Time\\(s)} & \shortstack{B\&B\\nds} & \shortstack{Time\\(s)} & \shortstack{B\&B\\nds} & \shortstack{Time\\(s)} & \shortstack{B\&B\\nds}  & \shortstack{STSP} & \shortstack{ATSP} \\
\noalign{\smallskip}\hline\noalign{\smallskip}
7  & 0.78   & 1    & 0.01   & 9     & 0.79   & 1    & 0.01  & 9    & 0.00    & 20     & 0.00    & 9     & 1.021 & 1.176 \\
8  & 3.33   & 1    & 0.02   & 11    & 3.45   & 1    & 0.12  & 11   & 0.01 & 17     & 0.01 & 9     & 1.014 & 1.154 \\
9  & 3.59   & 5    & 0.11   & 43    & 3.01   & 7    & 0.10   & 35   & 0.01 & 71     & 0.02 & 85    & 1.017 & 1.197 \\
10 & 5.35   & 13   & 0.13   & 61    & 3.22   & 11   & 0.11  & 57   & 0.01 & 126    & 0.02 & 73    & 1.017 & 1.174 \\
11 & 6.66   & 27   & 0.38   & 91    & 6.27   & 31   & 0.27  & 71   & 0.01 & 636    & 0.04 & 187   & 1.014 & 1.195 \\
12 & 5.62   & 41   & 0.95   & 289   & 11.72  & 117  & 0.69  & 137  & 0.01 & 1726   & 0.20  & 939   & 1.014 & 1.212 \\
13 & 25.33  & 217  & 1.18   & 135   & 14.52  & 165  & 0.83  & 125  & 0.03 & 5735   & 0.52 & 2271  & 1.015 & 1.248 \\
14 & 22.53  & 285  & 8.99   & 1465  & 21.34  & 265  & 4.75  & 812  & 0.07 & 15835  & 2.45 & 9923  & 1.013 & 1.258 \\
15 & 19.28  & 157  & 0.81   & 85    & 25.51  & 203  & 0.91  & 132  & 0.11 & 21238  & 0.81 & 2985  & 1.013 & 1.246 \\
16 & 58.62  & 575  & 15.01  & 1461  & 48.38  & 437  & 3.29  & 385  & 0.41 & 75912  & 8.36 & 28031 & 1.012 & 1.247 \\
17 & 314.79 & 3191 & 32.63  & 2323  & 90.09  & 1089 & 9.01  & 939  & 1.08 & 218251 & -    & -     & 1.012 & 1.258 \\
18 & 508.98 & 5737 & 559.31 & 26613 & 515.16 & 5999 & 97.84 & 8701 & 3.16 & 517070 & -    & -     & 1.010 & 1.227 \\
19 & 21.22  & 135  & 4.49   & 457   & 23.75  & 145  & 1.15  & 121  & 2.02 & 169364 & 3.62 & 9221  & 1.009 & 1.262\\
\noalign{\smallskip}\hline\noalign{\smallskip}
\end{tabular}
\end{adjustbox}
\end{minipage}
\end{center}
\end{sidewaystable}

\begin{sidewaystable}
%\sidewaystable
\begin{center}
\begin{minipage}{\textheight}
\footnotesize
\centering
\caption{Results of the runs of different solvers on the ATSPLIB \cite{reinelt1991tsplib}. The first two columns represent the number of nodes in the ATSP instances and the name as reported in ATSPLIB; \texttt{(-)} indicates that Concorde was run with the default settings; \texttt{-d} refers to running Concorde with the DF option enabled; \texttt{-C0} indicates that local cuts were disabled; \texttt{-d -C0} corresponds to running Concorde with both DF enabled and local cuts disabled. \texttt{ATSP92} reports the performance from \revv{\cite{fischetti1992additive}}, and \texttt{FLT07} shows the performance of the algorithm presented \revv{\cite{fischetti2007exact}}. The last column shows the integrality gap (\texttt{IG}) for both instances—the original ATSP and the transformed STSP, obtained through the procedure described in Section \ref{sec:hard_ATSPLIB}. Gaps are rounded to the third decimal place for readability.}\label{tab:comp_hard_ATSPLIB}
\begin{adjustbox}{scale=0.8}
\begin{tabular}{llrrrrrrrrrrrrrr}
\hline\noalign{\smallskip}
Name & $n$ & \multicolumn{2}{c}{ \rev{-} } & \multicolumn{2}{c}{\rev{\texttt{-d}}} & \multicolumn{2}{c}{\rev{\texttt{-C0}}} & \multicolumn{2}{c}{\rev{\texttt{-C0 -d}}}  & \multicolumn{2}{c}{\rev{\texttt{ATSP92}}} & \multicolumn{2}{c}{\rev{\texttt{FLT07}}} & \multicolumn{2}{c}{\rev{\texttt{IG}}}\\
& & \shortstack{Time\\(s)} & \shortstack{B\&B\\nds} & \shortstack{Time\\(s)} & \shortstack{B\&B\\nds} & \shortstack{Time\\(s)} & \shortstack{B\&B\\nds} & \shortstack{Time\\(s)} & \shortstack{B\&B\\nds} & \shortstack{Time\\(s)} & \shortstack{B\&B\\nds} & \shortstack{Time\\(s)} & \shortstack{B\&B\\nds} & \shortstack{STSP} & \shortstack{ATSP} \\
\noalign{\smallskip}\hline\noalign{\smallskip}
\texttt{br17}    & 17  & 0.01 & 1  & 0.01 & 1  & 0.01 & 1  & 0.01 & 1  & 0.01   & 22      & 0.00    & 1  & 1.000     & 1.000     \\
\texttt{ftv33}   & 34  & 0.01 & 1  & 0.02 & 1  & 0.01 & 1  & 0.01 & 1  & 0.01   & 164     & 0.00    & 1    & 1.000     & 1.010 \\
\texttt{ftv35}   & 36  & 0.34 & 5  & 0.08 & 11 & 0.38 & 5  & 0.06 & 9  & 0.03   & 1108    & 0.01 & 11   & 1.001 & 1.011 \\
\texttt{ftv38}   & 39  & 0.99 & 5  & 0.11 & 11 & 1.32 & 9  & 0.12 & 11 & 0.04   & 1403    & 0.02 & 17   & 1.001 & 1.010 \\
\texttt{p43}     & 43  & 2.04 & 17 & 0.39 & 33 & 0.75 & 11 & 0.22 & 31 & 25.32* & 174562* & 0.34 & 277  & 1.000 & 1.002 \\
\texttt{ftv44}   & 45  & 0.35 & 3  & 0.22 & 15 & 0.45 & 5  & 0.23 & 19 & 0.02   & 340     & 0.04 & 35   & 1.001 & 1.018 \\
\texttt{ry48p}   & 48  & 1.41 & 3  & 0.24 & 15 & 1.75 & 5  & 0.24 & 17 & 2.24   & 25910   & 0.14 & 89   & 1.000 & 1.016 \\
\texttt{ftv47}   & 48  & 0.43 & 5  & 0.22 & 9  & 0.49 & 5  & 0.20  & 11 & 0.14   & 5495    & 0.13 & 165  & 1.001 & 1.016 \\
\texttt{ft53}    & 53  & 0.03 & 1  & 0.02 & 1  & 0.02 & 1  & 0.03 & 1  & 0.02   & 141     & 0.01 & 1    & 1.000     & 1.000     \\
\texttt{ftv55}   & 56  & 1.06 & 5  & 0.18 & 5  & 1.00    & 5  & 0.18 & 5  & 0.21   & 9245    & 0.08 & 107  & 1.001 & 1.015 \\
\texttt{ftv64}   & 65  & 1.19 & 5  & 0.15 & 15 & 1.78 & 7  & 0.17 & 17 & 0.41   & 8224    & 0.12 & 119  & 1.001 & 1.018 \\
\texttt{ft70}    & 70  & 0.29 & 3  & 0.10  & 7  & 0.29 & 1  & 0.10  & 7  & 0.03   & 276     & 0.01 & 7    & 1.000 & 1.001 \\
\texttt{ftv70}   & 71  & 1.05 & 7  & 0.20  & 21 & 1.25 & 5  & 0.21 & 21 & 1.29   & 21606   & 0.35 & 277  & 1.001 & 1.021 \\
\texttt{kro124p} & 100 & 0.32 & 1  & 0.58 & 17 & 0.30  & 1  & 0.58 & 23 & 6.33   & 17387   & 0.18 & 59   & 1.000 & 1.009\\
\noalign{\smallskip}\hline\noalign{\smallskip}
\end{tabular}
\end{adjustbox}
\end{minipage}
\end{center}
\end{sidewaystable}

\color{black}
%%%%%%%%%%%%%%%%%%%%%%%%%%%%%%%%%%%%
\section{Conclusions}\label{sec:final}
In this paper, we have introduced and implemented a new symmetry-breaking pivoting algorithm and a new $\lambda$-loop breaking procedure that \rev{allow} the exploration of vertices of the asymmetric subtour elimination polytope, yielding a large integrality gap.
The symmetry-breaking pivoting exploits the class of \rev{isomorphisms} of the vertices of $\pasep{n}$ that we completely calculated for a small value of $n$.
Checking whether two vertices of $\pasep{n}$ are isomorphic is currently one of the two computational bottlenecks of our procedure.
For each non-isomorphic vertex we visit, we solve an instance of the Gap$(\bm x)$ problem.
With our new algorithm, we can compute new lower bounds for $\alpha^{LB}_n$ for $n\leq 22$ by exploring not only half-integer vertices. 

\rev{As a byproduct of our generation, we propose 13 small ATSP instances that are challenging for all the solvers we tested}

In the future, we plan to explore the unresolved issue addressed in this study, which involves developing a procedure that yields the stabilizer based on a vertex and creating a strategy that leverages symmetries to produce vertices that are considered \emph{noteworthy} in terms of the integrality gap.
Note that for small values of $n$, the integrality gap values returned by the two families of instances proposed and studied in \cite{charikar2004integrality,elliott2008integrality} are improved. For these two families, the integrality gap converges to $2$. Therefore, if the improvement obtained in this paper for small values of $n$ could be ``uniformly'' observed also for large values of $n$\rev{, this} would lead to an integrality gap greater than $2$. Of course, this is just an intriguing hypothesis for future research.

\rev{Another interesting research direction would be to study the integrality gap of formulations stronger than the DFJ formulation; that is, to analyze how the integrality gap evolves when additional valid inequalities, such as the Comb inequalities \cite{chvatal1973edmonds}, are added. In principle, our approach can be extended to work on other polytopes, and thus it could provide insights into the worst-case instances for the integrality gap under these stronger formulations.}

\section*{Acknowledgments}
The work of M. Mastrolilli, L. M. Gambardella, and E. Vercesi has been supported by the Swiss National Science Foundation project n. 200021\_212929/1 ``Computational methods for integrality gaps analysis'', Project code 36RAGAP. 

S. Gualandi acknowledges the contribution of the National Recovery and Resilience Plan, Mission 4 Component 2 - Investment 1.4 - NATIONAL CENTER FOR HPC, BIG DATA AND QUANTUM COMPUTING, spoke 6. 

We thank Giovanni Rinaldi for pointing out the interesting references \cite{padberg_rao} and \cite{balinski1974on}. 
\rev{We are deeply indebted to Matteo Fischetti for sharing the code from \cite{fischetti1992additive} and \cite{fischetti2007exact}, and for assisting us in getting it to compile and run on our machines.}

\newpage
\appendix
\section{Detailed description of the orbits for for $n = 6$}\label{appendix_n6}
\begin{table}[h]
\caption{Orbit structure for $n = 6$, top 21 having the highest integrality gap. Columns: cardinality of the orbit, type of components, frequency of each component, and integrality gap attained at the elements of that orbit.}\label{tab:orbit_6_1}
    \centering\medskip
    \begin{tabular}{rrrrrrrrrrrr}
    \hline\noalign{\smallskip}
        $\vert O_{\bm{x}} \vert $ & \multicolumn{5}{c}{Components} & \multicolumn{5}{c}{Frequencies} & IG   \\\\
        \noalign{\smallskip}\hline\noalign{\smallskip}
         180 & 0 & $1/2$ & ~ & ~ & ~ & 18 & 12 & ~ & ~ & ~ & $4/3$ \\
        120 & 0 & $1/3$ & $2/3$ & ~ & ~ & 18 & 6 & 6 & ~ & ~ & $9/7$ \\
        360 & 0 & $1/2$ & ~ & ~ & ~ & 18 & 12 & ~ & ~ & ~ & $14/11$ \\
        360 & 0 & $1/2$ & ~ & ~ & ~ & 18 & 12 & ~ & ~ & ~ & $5/4$ \\
        720 & 0 & $1/2$ & ~ & ~ & ~ & 18 & 12 & ~ & ~ & ~ & $5/4$ \\
        180 & 0 & $1/2$ & ~ & ~ & ~ & 18 & 12 & ~ & ~ & ~ & $5/4$ \\
        720 & 0 & 1 & $1/2$ & ~ & ~ & 19 & 1 & 10 & ~ & ~ & $5/4$ \\
        720 & 0 & 1 & $1/2$ & ~ & ~ & 19 & 1 & 10 & ~ & ~ & $5/4$ \\
        360 & 0 & 1 & $1/2$ & ~ & ~ & 19 & 1 & 10 & ~ & ~ & $5/4$ \\
        720 & 0 & $1/4$ & $3/4$ & $1/2$ & ~ & 17 & 5 & 3 & 5 & ~ & $16/13$ \\
        720 & 0 & $1/4$ & $3/4$ & $1/2$ & ~ & 16 & 6 & 2 & 6 & ~ & $6/5$ \\
        720 & 0 & $1/4$ & $3/4$ & $1/2$ & ~ & 16 & 6 & 2 & 6 & ~ & $6/5$ \\
        720 & 0 & $1/4$ & $3/4$ & $1/2$ & ~ & 17 & 5 & 3 & 5 & ~ & $6/5$ \\
        720 & 0 & $1/4$ & $3/4$ & $1/2$ & ~ & 16 & 6 & 2 & 6 & ~ & $6/5$ \\
        720 & 0 & $1/4$ & $3/4$ & $1/2$ & ~ & 16 & 6 & 2 & 6 & ~ & $6/5$ \\
        720 & 0 & $1/4$ & $3/4$ & $1/2$ & ~ & 16 & 6 & 2 & 6 & ~ & $6/5$ \\
        720 & 0 & $1/4$ & $3/4$ & $1/2$ & ~ & 16 & 6 & 2 & 6 & ~ & $6/5$ \\
        720 & 0 & $1/4$ & $3/4$ & $1/2$ & ~ & 16 & 6 & 2 & 6 & ~ & $6/5$ \\
        720 & 0 & $1/4$ & $3/4$ & $1/2$ & ~ & 16 & 6 & 2 & 6 & ~ & $6/5$ \\
        720 & 0 & $1/4$ & $3/4$ & $1/2$ & ~ & 16 & 6 & 2 & 6 & ~ & $6/5$ \\
        720 & 0 & $1/4$ & $3/4$ & $1/2$ & ~ & 17 & 5 & 3 & 5 & ~ & $6/5$ \\
        \noalign{\smallskip}\hline
    \end{tabular}
\end{table}

\newpage
\begin{table}[h]
\caption{Orbit structure for $n = 6$, bottom 69 having the lowest integrality gap. Columns: cardinality of the orbit, type of components, frequency of each component, and integrality gap attained at the elements of that orbit.}\label{tab:orbit_6_2}
    \tiny	
    \centering\medskip
    \begin{tabular}{rrrrrrrrrrrr}
    \hline\noalign{\smallskip}
        $\vert O_{\bm{x}} \vert $ & \multicolumn{5}{c}{Components} & \multicolumn{5}{c}{Frequencies} & IG   \\
        \noalign{\smallskip}\hline\noalign{\smallskip}
         720 & 0 & $1/3$ & $2/3$ & ~ & ~ & 17 & 8 & 5 & ~ & ~ & $6/5$ \\
        720 & 0 & $1/3$ & $2/3$ & ~ & ~ & 16 & 10 & 4 & ~ & ~ & $6/5$ \\
        720 & 0 & $1/3$ & $2/3$ & ~ & ~ & 16 & 10 & 4 & ~ & ~ & $6/5$ \\
        720 & 0 & $1/3$ & $2/3$ & ~ & ~ & 16 & 10 & 4 & ~ & ~ & $6/5$ \\
        720 & 0 & $1/3$ & $2/3$ & ~ & ~ & 16 & 10 & 4 & ~ & ~ & $6/5$ \\
        720 & 0 & $1/3$ & $2/3$ & ~ & ~ & 16 & 10 & 4 & ~ & ~ & $6/5$ \\
        720 & 0 & $1/3$ & $2/3$ & ~ & ~ & 16 & 10 & 4 & ~ & ~ & $6/5$ \\
        720 & 0 & $1/3$ & $2/3$ & ~ & ~ & 16 & 10 & 4 & ~ & ~ & $6/5$ \\
        720 & 0 & $1/3$ & $2/3$ & ~ & ~ & 16 & 10 & 4 & ~ & ~ & $6/5$ \\
        720 & 0 & $1/3$ & $2/3$ & ~ & ~ & 16 & 10 & 4 & ~ & ~ & $6/5$ \\
        720 & 0 & $1/3$ & $2/3$ & ~ & ~ & 16 & 10 & 4 & ~ & ~ & $6/5$ \\
        720 & 0 & $1/3$ & $2/3$ & ~ & ~ & 17 & 8 & 5 & ~ & ~ & $6/5$ \\
        720 & 0 & $1/3$ & $2/3$ & ~ & ~ & 17 & 8 & 5 & ~ & ~ & $6/5$ \\
        720 & 0 & $1/3$ & $2/3$ & ~ & ~ & 17 & 8 & 5 & ~ & ~ & $6/5$ \\
        720 & 0 & $1/3$ & $2/3$ & ~ & ~ & 17 & 8 & 5 & ~ & ~ & $6/5$ \\
        720 & 0 & $1/3$ & $2/3$ & ~ & ~ & 17 & 8 & 5 & ~ & ~ & $6/5$ \\
        720 & 0 & $1/3$ & $2/3$ & ~ & ~ & 17 & 8 & 5 & ~ & ~ & $6/5$ \\
        720 & 0 & $1/3$ & $2/3$ & ~ & ~ & 16 & 10 & 4 & ~ & ~ & $6/5$ \\
        720 & 0 & $1/3$ & $2/3$ & ~ & ~ & 17 & 8 & 5 & ~ & ~ & $6/5$ \\
        720 & 0 & $1/3$ & $2/3$ & ~ & ~ & 17 & 8 & 5 & ~ & ~ & $6/5$ \\
        720 & 0 & $1/3$ & $2/3$ & ~ & ~ & 17 & 8 & 5 & ~ & ~ & $6/5$ \\
        720 & 0 & $1/3$ & $2/3$ & ~ & ~ & 17 & 8 & 5 & ~ & ~ & $6/5$ \\
        720 & 0 & $1/3$ & $2/3$ & ~ & ~ & 16 & 10 & 4 & ~ & ~ & $6/5$ \\
        720 & 0 & $1/3$ & $2/3$ & ~ & ~ & 17 & 8 & 5 & ~ & ~ & $6/5$ \\
        720 & 0 & $1/3$ & $2/3$ & ~ & ~ & 17 & 8 & 5 & ~ & ~ & $6/5$ \\
        240 & 0 & $1/3$ & $2/3$ & ~ & ~ & 18 & 6 & 6 & ~ & ~ & $6/5$ \\
        720 & 0 & $1/3$ & $2/3$ & ~ & ~ & 16 & 10 & 4 & ~ & ~ & $6/5$ \\
        720 & 0 & $1/3$ & $2/3$ & ~ & ~ & 17 & 8 & 5 & ~ & ~ & $6/5$ \\
        720 & 0 & $1/3$ & $2/3$ & ~ & ~ & 17 & 8 & 5 & ~ & ~ & $6/5$ \\
        720 & 0 & $1/3$ & $2/3$ & ~ & ~ & 16 & 10 & 4 & ~ & ~ & $6/5$ \\
        120 & 0 & $1/2$ & ~ & ~ & ~ & 18 & 12 & ~ & ~ & ~ & $6/5$ \\
        720 & 0 & $1/2$ & ~ & ~ & ~ & 18 & 12 & ~ & ~ & ~ & $6/5$ \\
        720 & 0 & $1/2$ & ~ & ~ & ~ & 18 & 12 & ~ & ~ & ~ & $6/5$ \\
        720 & 0 & $1/2$ & ~ & ~ & ~ & 18 & 12 & ~ & ~ & ~ & $6/5$ \\
        720 & 0 & $1/2$ & ~ & ~ & ~ & 18 & 12 & ~ & ~ & ~ & $6/5$ \\
        720 & 0 & 1 & $1/3$ & $2/3$ & ~ & 18 & 1 & 7 & 4 & ~ & $6/5$ \\
        720 & 0 & 1 & $1/3$ & $2/3$ & ~ & 18 & 1 & 7 & 4 & ~ & $6/5$ \\
        720 & 0 & 1 & $1/3$ & $2/3$ & ~ & 18 & 1 & 7 & 4 & ~ & $6/5$ \\
        720 & 0 & 1 & $1/3$ & $2/3$ & ~ & 18 & 1 & 7 & 4 & ~ & $6/5$ \\
        360 & 0 & 1 & $1/2$ & ~ & ~ & 19 & 1 & 10 & ~ & ~ & $6/5$ \\
        720 & 0 & 1 & $1/2$ & ~ & ~ & 19 & 1 & 10 & ~ & ~ & $6/5$ \\
        720 & 0 & 1 & $1/2$ & ~ & ~ & 19 & 1 & 10 & ~ & ~ & $6/5$ \\
        720 & 0 & 1 & $1/3$ & $2/3$ & ~ & 18 & 1 & 7 & 4 & ~ & $6/5$ \\
        360 & 0 & 1 & $1/2$ & ~ & ~ & 20 & 2 & 8 & ~ & ~ & $6/5$ \\
        720 & 0 & 1 & $1/2$ & ~ & ~ & 20 & 2 & 8 & ~ & ~ & $6/5$ \\
        720 & 0 & 1 & $1/2$ & ~ & ~ & 20 & 2 & 8 & ~ & ~ & $6/5$ \\
        720 & 0 & $1/5$ & $4/5$ & $2/5$ & $3/5$ & 16 & 5 & 2 & 4 & 3 & $25/21$ \\
        720 & 0 & $1/5$ & $4/5$ & $2/5$ & $3/5$ & 16 & 5 & 2 & 4 & 3 & $25/21$ \\
        720 & 0 & $1/5$ & $4/5$ & $2/5$ & $3/5$ & 16 & 5 & 2 & 4 & 3 & $25/21$ \\
        720 & 0 & $1/4$ & $3/4$ & $1/2$ & ~ & 16 & 7 & 3 & 4 & ~ & $20/17$ \\
        720 & 0 & $1/4$ & $3/4$ & $1/2$ & ~ & 16 & 7 & 3 & 4 & ~ & $20/17$ \\
        720 & 0 & $1/4$ & $3/4$ & $1/2$ & ~ & 16 & 7 & 3 & 4 & ~ & $20/17$ \\
        720 & 0 & $1/4$ & $3/4$ & $1/2$ & ~ & 16 & 7 & 3 & 4 & ~ & $20/17$ \\
        720 & 0 & $1/4$ & $3/4$ & $1/2$ & ~ & 16 & 7 & 3 & 4 & ~ & $20/17$ \\
        720 & 0 & $1/4$ & $3/4$ & ~ & ~ & 16 & 9 & 5 & ~ & ~ & $20/17$ \\
        720 & 0 & $1/3$ & $2/3$ & ~ & ~ & 16 & 10 & 4 & ~ & ~ & $7/6$\\
        120 & 0 & $1/2$ & ~ & ~ & ~ & 18 & 12 & ~ & ~ & ~ & $7/6$ \\
        720 & 0 & $1/3$ & $2/3$ & ~ & ~ & 16 & 10 & 4 & ~ & ~ & $15/13$ \\
        720 & 0 & $1/3$ & $2/3$ & ~ & ~ & 17 & 8 & 5 & ~ & ~ & $15/13$ \\
        360 & 0 & $1/3$ & $2/3$ & ~ & ~ & 16 & 10 & 4 & ~ & ~ & $15/13$ \\
        720 & 0 & $1/3$ & $2/3$ & ~ & ~ & 16 & 10 & 4 & ~ & ~ & $15/13$ \\
        720 & 0 & $1/3$ & $2/3$ & ~ & ~ & 16 & 10 & 4 & ~ & ~ & $15/13$ \\
        720 & 0 & $1/3$ & $2/3$ & ~ & ~ & 16 & 10 & 4 & ~ & ~ & $15/13$ \\
        720 & 0 & $1/3$ & $2/3$ & ~ & ~ & 16 & 10 & 4 & ~ & ~ & $15/13$ \\
        360 & 0 & $1/3$ & $2/3$ & ~ & ~ & 16 & 10 & 4 & ~ & ~ & $15/13$ \\
        720 & 0 & $1/3$ & $2/3$ & ~ & ~ & 17 & 8 & 5 & ~ & ~ & $15/13$ \\
        720 & 0 & $1/3$ & $2/3$ & ~ & ~ & 17 & 8 & 5 & ~ & ~ & $15/13$ \\
        120 & 0 & $1/3$ & $2/3$ & ~ & ~ & 18 & 6 & 6 & ~ & ~ & $15/13$ \\
        120 & 0 & 1 & ~ & ~ & ~ & 24 & 6 & ~ & ~ & ~ & 1 \\
        \noalign{\smallskip}\hline
    \end{tabular}
\end{table}

%Bibliography
\printbibliography

@article{BoydE05,
  author       = {Sylvia C. Boyd and
                  Paul Elliott{-}Magwood},
  title        = {Computing the integrality gap of the asymmetric travelling salesman
                  problem},
  journal      = {Electron. Notes Discret. Math.},
  volume       = {19},
  pages        = {241--247},
  year         = {2005}
}

@phdthesis{elliott2008integrality,
  title={{The integrality gap of the Asymmetric Travelling Salesman Problem}},
  author={Elliott-Magwood, Paul},
  year={2008},
  school={University of Ottawa (Canada)}
}

@article{hougardy2020hard,
  title={{Hard to solve instances of the Euclidean Traveling Salesman Problem}},
  author={Hougardy, Stefan and Zhong, Xianghui},
  journal={Mathematical Programming Computation},
  volume={13},
  pages={51--74},
  year={2021},
  publisher={Springer}
}

@article{zhong2021lower,
title = {Lower bounds on the integrality ratio of the subtour {LP} for the traveling salesman problem},
	volume = {365},
	issn = {0166-218X},
	journal = {Discrete Applied Mathematics},
	author = {Zhong, Xianghui},
	year = {2025},
	pages = {109--129},
}

@book{vaziranibook,
  title={Approximation Algorithms},
  author={Vazirani, Vijay V.},
  volume={1},
  year={2001},
  publisher={Springer},
  address={Berlin}
}

@article{applegate1998solution,
  title={On the solution of traveling salesman problems},
    journal = {Documenta Mathematica},
  author={Applegate, David and Bixby, Robert and Cook, William and Chv{\'a}tal, Vasek},
  language = {eng},
pages = {645-656},
publisher = {Universiät Bielefeld, Fakultät für Mathematik},
title = {On the solution of traveling salesman problems.},
url = {http://eudml.org/doc/233207},
year = {1998},
}

@article{fischetti1997polyhedral,
    title={{A polyhedral approach to the Asymmetric Traveling Salesman Problem}},
    author={Fischetti, Matteo and Toth, Paolo},
    journal={Management Science},
    volume={43},
    number={11},
    pages={1520--1536},
    year={1997},
    publisher={INFORMS}
}

@article{fischetti2007exact,
  title={{Exact methods for the Asymmetric Traveling Salesman Problem}},
  author={Fischetti, Matteo and Lodi, Andrea and Toth, Paolo},
  journal={The traveling salesman problem and its variations},
  pages={169--205},
  year={2007},
  publisher={Springer}
}

@article{jonker1983transforming,
  title={{Transforming Asymmetric into Symmetric Traveling Salesman Problems}},
  author={Jonker, Roy and Volgenant, Ton},
  journal={Operations Research Letters},
  volume={2},
  number={4},
  pages={161--163},
  year={1983},
  publisher={Elsevier}
}

@article{roberti2012models,
  title={Models and algorithms for the {A}symmetric {T}raveling {S}alesman {P}roblem: an experimental comparison},
  author={Roberti, Roberto and Toth, Paolo},
  journal={EURO Journal on Transportation and Logistics},
  volume={1},
  number={1-2},
  pages={113--133},
  year={2012},
  publisher={Elsevier}
}

@article{benoit2008finding,
	author = {Benoit, Genevieve and Boyd, Sylvia},
	journal = {Mathematics of Operations Research},
	number = {4},
	pages = {921--931},
	publisher = {INFORMS},
	title = {Finding the exact integrality gap for small {T}raveling {S}alesman {P}roblems},
	volume = {33},
	year = {2008}}

@article{HOUGARDY2014495,
	author = {Stefan Hougardy},
	issn = {0167-6377},
	journal = {Operations Research Letters},
	number = {8},
	pages = {495--499},
	title = {{On the integrality ratio of the subtour LP for Euclidean TSP}},
	volume = {42},
	year = {2014},
}

@article{vercesi2023generation,
  title={{On the generation of metric TSP instances with a large integrality gap by branch-and-cut}},
  author={Vercesi, Eleonora and Gualandi, Stefano and Mastrolilli, Monaldo and Gambardella, Luca Maria},
  journal={Mathematical Programming Computation},
  volume={15},
  number={2},
  pages={389--416},
  year={2023},
  publisher={Springer}
}

@article{charikar2004integrality,
  title={On the integrality ratio for the {A}symmetric {T}raveling {S}alesman {P}roblem},
  author={Charikar, Moses and Goemans, Michel X and Karloff, Howard},
  journal={Mathematics of Operations Research},
  volume={31},
  number={2},
  pages={245--252},
  year={2006},
  publisher={INFORMS}
}

@article{carr2004held,
	title={{On the Held-Karp relaxation for the {A}symmetric and {S}ymmetric {T}raveling {S}alesman {P}roblems}},
	author={Carr, Robert and Vempala, Santosh},
	journal={Mathematical Programming},
	volume={100},
	number={3},
	pages={569--587},
	year={2004},
	publisher={Springer}
}

@article{fousse2007mpfr,
  title={{MPFR: A multiple-precision binary floating-point library with correct rounding}},
  author={Fousse, Laurent and Hanrot, Guillaume and Lef{\`e}vre, Vincent and P{\'e}lissier, Patrick and Zimmermann, Paul},
  journal={ACM Transactions on Mathematical Software (TOMS)},
  volume={33},
  number={2},
  pages={13--es},
  year={2007},
  publisher={ACM New York, NY, USA}
}

@MISC{eigenweb,
    author = {Ga\"{e}l Guennebaud and Beno\^{i}t Jacob},
    title = {Eigen v3},
    howpublished = {\url{http://eigen.tuxfamily.org}},
    year = {2010}
}

@misc{gurobi,
  author = {{Gurobi Optimization, LLC}},
  title = {{Gurobi Optimizer Reference Manual}},
  year = 2023,
  url = "https://www.gurobi.com"
}

@book{artin2011algebra,
  title={Algebra},
  author={Artin, Michael},
  year={2011},
  publisher={Pearson Education},
  address={Upper Saddle River, NJ}
}

@article{christof2009porta,
  title={Porta},
  author={Christof, Thomas},
  journal={http://www.iwr.uni-heidelberg. de/groups/comopt/software/PORTA/index. html},
  year={2009}
}

@article{lawvere1973metric,
  title = {Metric spaces, generalized logic, and closed categories},
  author = {Lawvere, F. William},
  journal = {Rendiconti del Seminario Matematico e Fisico di Milano},
  volume = {43},
  number = {1},
  pages = {135--166},
  year = {1973},
  month = {12},
}

@article{reinelt1991tsplib,
  title={{TSPLIB—A traveling salesman problem library}},
  author={Reinelt, Gerhard},
  journal={ORSA journal on computing},
  volume={3},
  number={4},
  pages={376--384},
  year={1991},
  publisher={INFORMS}
}

@inbook{Grotschelbook,
author    = "Grötschel, M. and Padberg, M. W.",
  editor = "Lawler, E. L. and Lenstra, J. K. and Rinnooy, A. H. G. and Shmoys, D. B.",
  title     = "The traveling salesman problem",
  chapter   = "Polyhedral theory",
  publisher = "John Wiley \& Sons Ltd",
  year      = "1985"
}

@article{padberg_rao,
  title = {The Travelling Salesman Problem and a Class of Polyhedra of Diameter Two},
  author = {Padberg, M.W. and Rao, M.R.},
  journal = {Mathematical Programming},
  volume = {7},
  pages = {32--45},
  year = {1974}
}

@article{singh2019integrality,
  title={Integrality gap of the vertex cover linear programming relaxation},
  author={Singh, Mohit},
  journal={Operations Research Letters},
  volume={47},
  number={4},
  pages={288--290},
  year={2019},
  publisher={Elsevier}
}

@article{assarf2017computing,
  title={Computing convex hulls and counting integer points with polymake},
  author={Assarf, Benjamin and Gawrilow, Ewgenij and Herr, Katrin and Joswig, Michael and Lorenz, Benjamin and Paffenholz, Andreas and Rehn, Thomas},
  journal={Math. Program. Comput.},
  volume={9},
  number={1},
  pages={1--38},
  year={2017},
  mr={3613012}
}

@article{mnich2018improved,
  title={Improved integrality gap upper bounds for traveling salesperson problems with distances one and two},
  author={Mnich, Matthias and M{\"o}mke, Tobias},
  journal={European journal of operational research},
  volume={266},
  number={2},
  pages={436--457},
  year={2018},
  publisher={Elsevier}
}

@inproceedings{boyd2011tsp,
  title={TSP on cubic and subcubic graphs},
  author={Boyd, Sylvia and Sitters, Ren{\'e} and van der Ster, Suzanne and Stougie, Leen},
  booktitle={Integer Programming and Combinatoral Optimization: 15th International Conference, IPCO 2011, New York, NY, USA, June 15-17, 2011. Proceedings 15},
  pages={65--77},
  year={2011},
  organization={Springer}
}

@article{dantzig1954solution,
  title={Solution of a large-scale traveling-salesman problem},
  author={Dantzig, George and Fulkerson, Ray and Johnson, Selmer},
  journal={Journal of the operations research society of America},
  volume={2},
  number={4},
  pages={393--410},
  year={1954},
  publisher={INFORMS}
}

@article{balinski1974on,
author = {Balinski, M. L. and Russakoff, Andrew},
title = {On the Assignment Polytope},
journal = {SIAM Review},
volume = {16},
number = {4},
pages = {516-525},
year = {1974},
doi = {10.1137/1016083},

URL = { 
    
        https://doi.org/10.1137/1016083
    
    

},
eprint = { 
    
        https://doi.org/10.1137/1016083
    
    

}
}

@inproceedings{traub2020improved,
  title={An improved approximation algorithm for ATSP},
  author={Traub, Vera and Vygen, Jens},
  booktitle={Proceedings of the 52nd annual ACM SIGACT symposium on theory of computing},
  pages={1--13},
  year={2020}
}

@article{fischetti1992additive,
  title={An additive bounding procedure for the asymmetric travelling salesman problem},
  author={Fischetti, Matteo and Toth, Paolo},
  journal={Mathematical programming},
  volume={53},
  pages={173--197},
  year={1992},
  publisher={Springer}
}

@manual{cplex,
  title        = {{IBM ILOG CPLEX Optimization Studio}},
  author       = {{IBM Corporation}},
  year         = {2019},
  note         = {Version 12.9.0},
  url          = {https://www.ibm.com/products/ilog-cplex-optimization-studio}
}

@article{chvatal1973edmonds,
  title={Edmonds polytopes and weakly Hamiltonian graphs},
  author={Chv{\'a}tal, V{\'a}clav},
  journal={Mathematical programming},
  volume={5},
  pages={29--40},
  year={1973},
  publisher={Springer}
}

\end{document}